\date {May 30, 2000} 
\title {  
Analyticity of intersection exponents \\ for planar Brownian motion
}
\author {Gregory F.  Lawler\thanks{Duke University, Research
supported by the National Science Foundation}
\and Oded Schramm\thanks{Microsoft Research}
\and Wendelin Werner\thanks{Universit\'e Paris-Sud}
}

\documentclass[12pt]{article}
\usepackage{amsmath}
\usepackage{amsthm}
\usepackage{amsfonts}
\usepackage{graphicx}
  \ifx\LabelFigloaded\MYundefined\relax
  \else
   \fi

  \def\LabelFigloaded{\relax}

  \chardef\LabelFigCatAt\the\catcode`\@
  \catcode`\@=11

 \let\LabelFigwlog@ld\wlog
 \def\wlog#1{\relax}

 \ifx\\\MYundefined@
    \let\\\relax
 \fi

  \def\ms@g{\immediate\write16}

 \def\N@wif{\csname newif\endcsname }
 \def\Temp@ {\N@wif\ifIN@}
 \ifx\INN@\MYundefined@
    \else \let\Temp@\relax
 \fi
 \Temp@

  \def\IN@{\expandafter\INN@\expandafter}
  \long\def\INN@0#1@#2@{\long\def\NI@##1#1##2##3\ENDNI@
    {\ifx\m@rker##2\IN@false\else\IN@true\fi}%
     \expandafter\NI@#2@@#1\m@rker\ENDNI@}
  \def\m@rker{\m@@rker}
 
  \newtoks\Initialtoks@  \newtoks\Terminaltoks@
  \def\SPLIT@{\expandafter\SPLITT@\expandafter}
  \def\SPLITT@0#1@#2@{\def\TTILPS@##1#1##2@{%
     \Initialtoks@{##1}\Terminaltoks@{##2}}\expandafter\TTILPS@#2@}

 \def\Shifted@@#1#2#3{\setbox0=\hbox{#3}%
   \raise -\dp0\vbox {\kern-#2%
       \hbox {\kern#1\unhbox0\kern-#1}%
           \kern#2}}

 \newcount\gridcount
 \newbox\auxGridbox@ \newbox\hGridbox@ \newbox\vGridbox@
 \newbox\Labelbox@ \newbox\auxLabelbox@
 \newbox\Coordinatebox@
 \newtoks\Labeltoks@
 \newdimen\Wdd@ \newdimen\Htt@
 \newdimen\Wddd@ \newdimen\Httt@
 
 \def\Wr@{\immediate\write16}

 \newdimen\GL@wd
 \def\GridLineWidth#1{\GL@wd=#1}

 \def\gobble#1{}
 \def\EdgeErr@{\Wr@{}%
      \Wr@{\string\Edges\space argument
      1, 10, 100 or 1000 please\string!}%
      }

 \newcount\Edgect@

 \def\Sweepup#1\endSweepup{}

 \def\SetEdges@{%
    \edef\Zr@@s{\expandafter\gobble\number\Edgect@\empty}%
        \count255=0\Zr@@s\relax
        \ifnum\count255=\z@\else\EdgeErr@\show\tailtest\fi
        \count255=1\Zr@@s\relax
        \ifnum\count255=\Edgect@\relax\else\EdgeErr@\show\leadtest\fi
    \EdgGl@b\edef\Zr@s{\expandafter\gobble\Zr@@s\empty}
    \ifnum\Edgect@>\@ne\relax\EdgGl@b\let\L@Dc\empty
        \else\EdgGl@b\edef\L@Dc{\string.}\fi
    \ifnum\Edgect@>\@ne\relax
        \EdgGl@b\edef\Edgescale@##1{\divide##1 by \Edgect@}%
        \else\EdgGl@b\edef\Edgescale@##1{}\fi
    }

 \def\Edges#1{\Edgect@=#1\relax
     \let\EdgGl@b\global \SetEdges@}

 \Edges{1}

 \def\hhrule{\hrule height \GL@wd\vskip-.\GL@wd}

 \def\hRule@{%
   \advance\gridcount -2%
   \vfil\hhrule\vfil
   \llap{\smash{\raise -2.5pt
     \hbox{\L@Dc\number\gridcount\Zr@s\kern2pt}}}%
   \hhrule
   }

\def\vvrule{\vrule width \GL@wd \kern-\GL@wd}

 \def\vRule@{\advance\gridcount 2%
   \hfil\vvrule\hfil
   \setbox\auxGridbox@=\vbox to 0pt
      {\vskip \Htt@\vskip 2pt
        \hbox to 0pt{\hss\L@Dc\number\gridcount\Zr@s\hss}\vss}%
      \wd\auxGridbox@=0pt \box\auxGridbox@
   \vvrule
   }

 \def\PlaceGrid@@{\gridcount=10 
  \setbox\hGridbox@=\hbox{%
        \hbox{%
             \hskip-.4pt\vrule
             \vbox to \Htt@{%
               \offinterlineskip\parindent=\z@\relax
               \hbox to \Wdd@{\hfil}
               \hRule@\hRule@\hRule@\hRule@
               \vfil\hhrule\vfil}%
             \vrule\hskip-.4pt}
    }%
  \gridcount=0%
  \setbox\vGridbox@=\hbox{%
      \vbox{\offinterlineskip\parindent=0pt\hsize=0pt
         \vskip-.4pt\hrule%
         \hbox to \Wdd@{%
                 \vtop to \Htt@{\vfil}%
                 \vRule@\vRule@\vRule@\vRule@
                 \hfil\vvrule\hfil}%
         \hrule\vskip-.4pt}}%
  \wd\hGridbox@=0pt\ht\hGridbox@=0pt
  \wd\vGridbox@=0pt\ht\vGridbox@=0pt
  \hbox{\box\hGridbox@\box\vGridbox@}%
  }

 \def\LabelsGlobal{\def\LabGl@b{\global}}
 \def\LabelsLocal{\def\LabGl@b{}}
 \LabelsGlobal 

 \def\SetLabels#1\endSetLabels{%
   \LabGl@b\Labeltoks@={#1()\\}%
   }

 \LabGl@b\Labeltoks@={()\\}

 \def\ShowGrid{\LabGl@b\let\PlaceGrid@\PlaceGrid@@}
 \def\HideGrid{\LabGl@b\let\PlaceGrid@\relax}
 \def\Grids{\ShowGrid\LabGl@b\let\GridSwitch@\ShowGrid}
 \def\noGrids{\HideGrid\LabGl@b\let\GridSwitch@\HideGrid}

 \noGrids

 \def\bAdjust@@{%
     \setbox\auxLabelbox@=\hbox{\raise \dp\auxLabelbox@
            \box\auxLabelbox@}}
 \def\bAdjust@{\let\vAdjust@\bAdjust@@}

 \def\eAdjust@@{\dimen0=-.5\ht\auxLabelbox@
     \advance\dimen0 by .5\dp\auxLabelbox@
     \setbox\auxLabelbox@=
            \hbox{\raise\dimen0\box\auxLabelbox@}}
 \def\eAdjust@{\let\vAdjust@\eAdjust@@}

 \def\tAdjust@@{%
     \setbox\auxLabelbox@=\hbox{\raise-\ht\auxLabelbox@
            \box\auxLabelbox@}}
 \def\tAdjust@{\let\vAdjust@\tAdjust@@}

 \let\vAdjust@\relax

 \def\lAdjust@{\let\hAdjust@\rlap}
 \def\rAdjust@{\let\hAdjust@\llap}

 \let\hAdjust@\relax\let\vAdjust@\relax

 \def\FetchLabel@#1(#2)#3\\{%
     \IN@0#2@@\ifIN@
        \setbox0=\hbox{\ignorespaces#1#3\unskip}%
        \ifdim\wd0>0pt
           \ms@g{}%
           \ms@g{ !!! Bad label(s)? !!!}%
           \message{ #1(#2)#3}%
        \fi
        \def\LabelMole@##1\endFetchLabel@{%
            \IN@0()\\@##1@%
            \ifIN@\def\Temp@{\FetchLabel@##1\endFetchLabel@}%
            \else\def\Temp@{}%
            \fi
            \Temp@
           }%
     \else
       \ignorespaces#1\unskip
       \setbox\auxLabelbox@=%
         \hbox to 0pt{\hss\ignorespaces\hAdjust@
          {\ignorespaces#3\unskip}\hss}%
       \vAdjust@
       \let\hAdjust@\relax\let\vAdjust@\relax
       \AugmentLabelBox@@{#2}%
       \ht\Labelbox@=0pt\dp\Labelbox@=0pt
       \let\LabelMole@\FetchLabel@%
     \fi\LabelMole@}

 \newtoks\XYSep@ 
 \def\SetXYSeparator#1{%
     \IN@0#1@@\ifIN@\XYSep@{*}%
     \else
     \XYSep@{#1}%
     \fi
     }

 \SetXYSeparator*

 \def\AugmentLabelBox@@#1{%
     \IN@0\the\XYSep@ @#1@\ifIN@
       \SPLIT@0\the\XYSep@ @#1@%
       \setbox\Labelbox@=\hbox to 0pt{%
         \unhbox\Labelbox@
         \Shifted@@{\the\Initialtoks@\Wddd@}%
         {\the\Terminaltoks@\Httt@}%
         {\box\auxLabelbox@}}%
     \else
         \ms@g{}%
         \ms@g{ !!! Bad insertion point. !!!}%
         \message{ (#1\ this point was rejected.)}%
     \fi
    }

 \def\FetchOption@#1[#2]#3\endFetchOption@{%
    \def\temp{#1}
    \ifx\temp\empty
       \Edgect@=#2\relax
       \let\EdgGl@b\relax
       \SetEdges@
       \Cleaner@#3%
    \fi}

 \def\Cleaner@#1[@]{\Labeltoks@{#1}}
     
 \def\PlaceLabels@@{\mathsurround=0pt
     \def\Cr@{\\}%
     \let\L\lAdjust@\let\R\rAdjust@
     \let\B\bAdjust@\let\E\eAdjust@\let\T\tAdjust@
     \expandafter\FetchOption@\the\Labeltoks@[@]\endFetchOption@
     \Wddd@=\Wdd@ \Edgescale@\Wddd@ 
     \Httt@=\Htt@ \Edgescale@\Httt@
     \expandafter\FetchLabel@\the\Labeltoks@\endFetchLabel@
     \box\Labelbox@
     }%

 \let \PlaceLabels@\PlaceLabels@@

 \def\AffixLabels#1{\setbox\Coordinatebox@=\hbox{#1}%
      \Wdd@=\wd\Coordinatebox@ \Htt@=\ht\Coordinatebox@
      \advance\Htt@ \dp\Coordinatebox@
      \hbox{\copy\Coordinatebox@\kern-\Wdd@ 
           \Shifted@@{0pt}{-\dp\Coordinatebox@}%
           {\PlaceLabels@\PlaceGrid@}%
           \kern\Wdd@}%
      \GridSwitch@ 
      \LabGl@b\Labeltoks@{()\\}%
      }
 
   \let\wlog\LabelFigwlog@ld   
   \catcode`\@=\LabelFigCatAt  

                                By

              Raymond S\'eroul <A18645@FRCCSC21.BITNET>
                                and 
              Laurent Siebenmann <lcs@topo.math.u-psud.fr>
    
              VERSIONS: July 1991, Oct 1991, Jan 1992, July 1992

INTRODUCTION

      This labelling package is intended for TeX users who
rely on non-TeX sources for for their graphics inserts.  It
provides means for adding TeX labels to such inserts with a
minimum of fuss. 

       For most labels, TeX users have in the past found it
reasonably convenient to rely on non-TeX sources. Typical
occasions when an inescapable need for TeX labels seemed to
arise are

 (a) when the graphics program lacks certain exotic or complex
mathematical symbols

 (b) when the very highest typographical quality is wanted for the
labels

 (c) when labels included with the graphics fail to print, 
 and you cannot figure out why (cf. boxedeps.doc).  The labels
 provided by labelfig.tex are 100

       Since this package first appeared, many users, who in the
past scarcely dreamed of using TeX labels, have come to use
nothing but.  So it is now appropriate to add

Intoxication Warning:  TeX labels may be addictive and expensive. 

     If you have a fast preview you may disagree, and even find
that this package provides an agreeable paste-up environment; see
extra applications at end.

     Note to publishers: It is possible and convenient to ultimately
export the TeX labels produced by labelfig.tex to become an integral
part of the EPS file. This is often desired by a publisher who typically
uses an "upmarket" graphics or page layout program, with which the
staff is skilled in perfecting figures.  See Appendix I for
a recipe.

     The authors are grateful to Patrick Ion of Math Reviews for
helpful comments and encouragement.

BASIC INSTRUCTIONS

    After reading in the macro file using

preview or proof your figure with a coordinate grid printed on
top, by typing the following:

    \ShowGrid  
    \AffixLabels{<the graphics insertion>}

Here <the graphics insertion> is what you would type to insert
the graphics object alone without the grid.  This must provide
for the space around it. For example <the graphics insertion>
might well be \BoxedEPSF{MyFigure scaled 700} using the
boxedeps.tex macro package (from same source); this provides a
TeX box containing the encapsulated PostScript insert specified by
the file MyFigure. \AffixLabels{...} provides the grid (supposing
\ShowGrid is present) and later, once you have specified labels
using the grid, it will "tack on" the labels.

     The grid is a sort of (usually elongated) checkerboard of
ten rows and ten columns and its (internal) partitions are by
default numbered  .1, ... ,.9  both horizontally (X-coordinate
running left to right) and vertically (Y-coordinate running bottom
to top).  Thus the points enclosed by the grid correspond to the
points of the unit square in the cartesian "X-Y" plane, the lower
left corner corresponding to the origin (0,0).  By extrapolation,
the full page corresponds to a larger rectangle in the plane.

     These coordinates serve to position labels as follows.
Before the \AffixLabels{...} command type label specifications:

  \SetLabels
   (<X-coordinate>*<Y-coordinate>) <first label> \\
   .
   .
   .
   (<X-coordinate>*<Y-coordinate>)  <last label> \\
  \endSetLabels

Each row specifies one label and is terminated by \\.  In each
row, the position indicator comes first; it is written as a
standard cartesian point except that the X- and Y- coordinates
are separated by * rather than a comma because TeX allows a
comma as decimal point. There are no dimension units to specify
as the unit is the grid itself.

     By default, this cartesian point specifies where the middle
of the baseline of the label will be located.  However if you precede
the point by \L [or \R] the left [or right] edge of the baseline will
be located there. Similarly you may also precede the point by \T, \E,
or \B to vertically align the top equator or bottom of the label box
at the specified point.  This gives nine standard positions of
the label with respect to the insertion point --- corresponding to
the eight principle points of the compas and the center

                     \L\T     \T      \R\T

                     \L\E     \E      \R\E

                     \L\B     \B      \R\B

But this neglects the default "baseline" level of TeX,
giving potentially three more positions

                     \L    <no tag>   \R

For text, the baseline level is often the preferred. Its relation to
the others is variable. It will often coincide with the bottom level,
as happens for "X".  But it is often distinct, as for "g", in which
case you have in all 12 distinct positions rather than 9.

     It is convenient to think of this specification of label
position as attaching the label by a thumb-tack to the coordinate
grid. There are up to twelve positions of the thumb-tack on the
label, while the position of the thumb-tack on the coordinate grid is
arbitrary.  Normally, one choses the position of the thumb-tack on
the label to be the one that is the closest to the item being
labeled.  There are good reasons for this "rule of thumb":

   (a)  It facilitates correct positioning at first try.

   (b)  If the scale of the figure must be altered after labels
have been affixed, the labels have a good chance of remaining well
positioned.

   (c)  The visible grid need not extend beyond the "bounding box"
for the figure, because the best preferred position is always
(at least almost) within the bounding box .

The second reason is particularly important. Indeed it often
happens that scale has to be altered after labelling begins, in
order to either provide space for the labels, or to adjust
proportions between the labels and the figure.  (The size of labels
is unaffected by scaling.)

     Here is an artificial but self-contained test which uses
TeX rules to make a graphics object.

TEST

    Do not skip this!

 \def\FrameIt#1{\hbox{\vrule$\vcenter {\hrule\kern3pt%
             \hbox {\kern3pt #1\kern3pt}%
               \kern3pt\hrule}$\relax\vrule}}

 \def\Caption#1#2{\FrameIt{%
       \vtop {\hsize=#1\relax \parindent=0pt
         \leftskip=0pt \rightskip=0pt plus15pt
         \parfillskip=0pt
         \lineskip=1pt\baselineskip=0pt
         #2}}}

 \def\FirstQuadrant{\hbox to 100pt{\vrule\vbox to 100pt{%
        \hbox to 100pt{\hfil}\vfil\hrule}\hss}}

  \SetLabels
    \R(.5*.2) $\zeta\,\cdot$\\
    (.9*-.10) $\xi$\\
    \R(-.03*.9) $\eta$\\
    \T(.5*.9) \Caption{70pt}{%
          \it The norm of
          $g(\xi+i\eta)$ is indicated on
          contours of this invisible surface.}\\
  \endSetLabels

  \AffixLabels{\FirstQuadrant} \end

  Note that the coordinates to use for labels are indicated on the
edges of the grid (when visible) corresponding to the conventional
x- and y- axes of the Cartesian plane. By default the grid is
1-by-1. However, by the command \Edges{100}, you can change this
to 100-by-100 and many users find this alternative most
convenient. Place the command \Edges{...} in your style file (or
header) since its effect is is global. Other possible edge values
are 10 and 1000.

  If you use the command \Edges{...} at all, do so with care.  For
if you accidentally delete an \Edges{...} command your labels will
abruptly be badly misplaced and may logically but mysteriously
generate "dimension too big" errors under TeX and "off page" errors
under your driver.  

  You can dictate the edgescale for an individual figure by giving
the scale in brackets immediately after \SetLabels.  Thus, to
import into an article using say \Edge{100} a figure labelled using
another edgescale, say the original 1-by-1 default, you can use
\SetLabels[1]...\endSetLabels.

GETTING IT DOWN PAT

     Complicated labeling deserves the same respect as
complicated mathematics.  Do not expect it to come out perfect the
first time!  What is needed in either case is a mechanism to
repeatedly typeset troublesome pieces.

     One mechanism is always available.  One does complicated
labelling in a separate "test" file involving just the figure being
labelled;  a texpert will know how to \dump TeX's current state as
a temporary format that restarts rapidly at each retry.  Usually,
one then pastes the completed labelled figure back into the main
TeX file, but, of course, one can also \input it as an auxiliary
file.

     If you do not have a TeXpert at handy, here is a first
approximation to an efficient setup. By deletions reduce a copy
of your article to just a few lines before and after the figure.
Now label the figure, and finally, copy and paste the labelled
figure to the original article. Then copy the next figure to label
into this testbed and repeat. The TeXpert can improve the  speed
at which TeX starts up, by compiling a format specifically for
your article; just one caution: best NOT include in the format
ephemeral details of setup like \Set<mydriver>ArtSpecials (from
boxedeps.tex because this reads  figure dimensions which you may
change during your work session.

     An improved mechanism to repeatedly typeset troublesome
pieces is now available on the Macintosh; it is called LinoTeX;
see the same ftp sources.  It could be set up on many types
of computer.

     Before using labelfig.tex to attach labels to a graphics
object inserted using boxedeps.tex or BoxedArt.tex, make it a
firm rule to carefully adjust the bounding box using the trimming
commands of these packages, and also at least tentatively scale
and position the object. Beware of changing the grid inadvertently
after the labels have been positioned.  For example, correcting
the bounding box of a PostScript graphics object can foul up the
labels by changing the coordinate grid to which the labels are
attached. This is particularly true for the trimming  commands of
boxedeps.tex and BoxedArt.tex. However, as noted already, change
of scale is much less disruptive, and modest adjustments should be
well tolerated.

     Sometimes the labels protrude so far from the bounding box
of a figure that the figure has to be repositioned.  Best do this
by ad hoc spacing, say using \hglue and \vglue; altering the
bounding box would create a vicious circle.

     Remember that you are responsible for preventing labels
from overlapping. You are responsible for all label typography
including size and style. A label is really just about anything
that can be put in a TeX box. Note that spaces at the beginning
and end of labels will normally be suppressed; if you really want
them you must protect them with TeX braces.

     This package temporarily sets the \mathsurround parameter
of TeX to zero  while the labels are being affixed. This is done
because nonzero \mathsurround space would influence the position
of left and right aligned labels; then, when a texpert or printer
modifies mathsurround, diagram labeling might be disastrously
altered. There is a small price to pay involving labels that are
formatted as caption boxes including mathematics: you  may want or
need to specify an explicit mathsurround space within the caption
box; it will not influence anything outside.

     Those hostile to the use of * as separator between
the X and Y coordinates of label insertion points, are free to
impose another using \SetXYSeparator{<the new separator>}.  
Americans may prefer "," to "*" since they never use a 
comma as a decimal point; on the other hand, * may be more visible.

APPENDIX (I)  MERGING labelfig.tex LABELS INTO AN EPSF GRAPHICS OBJECT.

     As promised in the introduction, here is a recipe useful for
publishers. It works at least on Macintosh and at least for vectorized
graphics and Adobe type1 fonts.  (There is surely a similar recipe for
PCs under MSWindows.)

 (a)  Use boxedeps.tex utility to integrate the figure given by the eps
file, "x.eps" say, with a visible frame around it.  See
\ShowDisplacementBoxes command in boxedeps.tex.  To get precise results
automatically it is important to use the \Trim... commands of
boxedeps.tex making the "DisplacementBox" neatly fit the figure.

 (b)  Use the TeX printer driver and LaserWriter (versions >= 8.1.1) to
export to an EPSF the DVI page containing the integrated, labelled
figure. You now have an EPS file  "xx.eps"  that contains too much, and at
the wrong scale, and at wrong position.

 (c)  Convert the EPSF to an Adode Illustrator format EPSF using
the shareware utility called epsConvert by Sam Weiss
1993-- (currently $25).

 (d)  In Illustrator (or a compatible program), group the labels and the
"DisplacementBox"; copy them to the clipboard and paste them into "x.ps".
This step requires that all the label fonts be "visible to the Macintosh.

 (e)  Translate and scale the pasted group consisting of the labels plus
the "DisplacementBox" so as to make the "DisplacementBox" the bounding
box of (labelless) figure represented by "x.eps".  At this point the
labels will be correctly placed on the figure "x.eps".

 (f)  Ungroup and delete the "DisplacementBox".  The result is the
desired single EPS file, "x+.eps" say, It contains the original figure
plus its labels.  

     Using grouping and ungrouping appropriately in "x+.eps", a
publisher's staff can very efficiently improve label positions etc.

APPENDIX II)  SOME EXOTIC APPLICATIONS

     The grid of labelfig.tex is analogous to a light-table in
classical page makeup with wax or latex glue.  In principle, you
can use it to compose any page from its indivisible parts.  This
even has some of the artisanal charm of classical paste-up
provided you have a fast screen preview to make the process
"interactive".

     In practice labelfig.tex is a tool for nonstandard jobs.
Here are a few going beyond the labelling already discussed.

(I)  GRAPHICS INTEGRATION.

     This is accomplished by treating the imported graphics
objects as labels.  The underlying graphics object is then
typically an empty  \vbox to <dimension>{\vfill} in a TeX
\midinsert...\endinsert construction.  A label line
might be of the form

   (.1*.1) \\

The exact form of the special command varies from driver to
driver.  However, in the case of encapsulated PostScript graphics
(EPSF norm), by relying on boxedeps.tex, one can have the
following standard syntax (independant of driver  (see
boxedeps.doc for details.
  
  (.1*.1) \BoxedEPSF{MyFigure scaled <scale in mils>}\\

This may be slow since it requires TeX to read the PostScript
file to read bounding box using many complex macros.  So you
may want to try

  (.1*.1) \EPSFSpecial{MyFigure}{<scale in mils>}\\

which is fast and driver independant, but it squashes the
bounding box, normally to its lower left corner.

     Similarly for graphics of the Macintosh PICT norm ---
using BoxedArt.tex (same sources) in place of boxedeps.tex.

     This approach to integration is to be recommended when
one is assembling a composite graphics object.

 (II)  COMMUTATIVE DIAGRAM ENHANCEMENT

     Commutative diagrams or arrays of mathematical objects
connected by arrows of various sorts are common in mathematics.
The mathematical objects require the use of TeX.  Recently TeX
acquired a good collection of arrows of all slopes --- that of
LamSTeX --- plus pwerful macros to build the diagrams.

     However, even the LamSTeX collection is often
inadequate; it lacks for example double shafted arrows, dotted
arrows and curved arrows. Fortunately it is possible to produce
such arrows on an individual basis using sophisticated graphics
programs such as Illustrator and AldusFreehand (both serving
the EPSF norm) or using Metafont (with its public domain norm).
Since the creation of each new arrow is a work of love, you
probably want to limit the number of arrows by using LamSTeX
for most arrows. The 40K commutative diagram module of LamSTeX
has been adapted to work with AmSTeX and a copy may be posted
with LabelFig and related files. Unfortunately no one has yet
offered a version that works with Plain TeX or LaTeX.

       Suffice it here to say that when the exotic arrow has
been somehow imported into TeX, labelfig.tex treats it as a
label that one affixes to the commutative diagram.  Two other
steps will be treated in separate notes, namely the matter of
extracting the dimension specifications for the arrow and the
construction of the arrow --- for these steps are far from
unique and often depend intimately on your computer environment. 
Notes for the Macintosh-Textures-Illustrator combination are
found in the file ExoticArrows.doc.

 (III) NESTING 

Ingenuity pays off in exploiting labelfig.tex. One can
mix graphics and typography quite freely.  labelfig.tex is good
for freeform or overlapping arrangements, while boxedeps.tex (or
BoxedArt.tex) is best for regimented non-overlapping
arrangements --- and the two can be combined.

     The default behavior of labelfig.tex is not ideal 
for nesting objects, because to prevent trouble for beginners
the register for labels is globally cleared when \AffixLabels
concludes.  But there are switches available

      \LabelsGlobal      \LabelsLocal

which change this.  To understand this, extend the above test 
by something like:

 \LabelsLocal

 \SetLabels
    (.5*.5) AAA\\
 \endSetLabels

 {
 \SetLabels
    (.5*.5) ZZZ\\
 \endSetLabels
   \AffixLabels{\FirstQuadrant}
 }

   \AffixLabels{\FirstQuadrant}

     There are however potential pitfalls.  Neither
labelfig.tex nor boxedeps.tex has been tested under extreme
conditions. Problems may occur if their procedures are
indiscriminately nested. For boxedeps.tex (not labelfig.tex)
there is a precise cause for worry, namely many of its
variables are "global", which means that TeX braces will not
provide the protection one might expect.

COMMAND SUMMARY FOR labelfig.tex

  Here [...] means optional (one or zero)
       [...]* means any number of such constructs

  \SetLabels
    [[<P>](<X><Sep><Y>) <label> \\]*
  \endSetLabels
  \ShowGrid  
  \AffixLabels{<the figure>}

   --- <P> is tack position, one of eleven or empty
              order irrelevant

                   \L\T      \T      \R\T

                   \L\E      \E      \R\E

                     \L               \R

                   \L\B      \B      \R\B

   --- (<X><Sep><Y>) insertion point;
  <Sep> is separator, = * by default;
  \SetXYSeparator{<Sep>} changes it.
   <X> and <Y> are real numbers

  --- <label> a label to attach 

  --- <the figure> the figure to label 

  \GlobalLabels (default)     
  \LocalLabels  setting for nested constructs.

 \Grids makes ALL grids appear; \HideGrid then makes just next disappear.
 \noGrids returns to default.  The commands are always global.

 \GridLineWidth{<dimension>} adjusts width of grid lines. Default is very
small, to give "hairline" effect. If your grid lines are missing try
setting \GridLineWidth{1pt}.

 \Edges#1 globally changes the edge size of all grids to the numerical 
value #1, which must be 1, 10, 100, or 1000.  The default is 1.

VERSION HISTORY.
 --- Jan 1993: \Edges#1 and [??] option after \SetLabels
 --- July 1992: \Grids, \noGrids, \HideGrid;
       Gridlines become hairlines; \GridLineWidth{<dimension>}.
 --- Oct 1991, Jan 1992: \SetXYSeparator{<Sep>},  \LabelsGlobal,
       \LabelsLocal.
 --- July 1991: first release

Address for bugs and other feedback:

        Raymond S\'eroul
        IREM and Lab. de Typographie Informatise
        Univ. Rene Descartes
        Strasbourg

    Tel 33-88-41-63-45
    Email:  A18645@FRCCSC21.BITNET

        Laurent Siebenmann
        Mathematique, Bat. 425,
        Univ de Paris-Sud,
        91405-Orsay,
        France

    Tel 33-1-6941-7949; 
    Email: lcs@topo.math.u-psud.fr  

\newif\ifhyper\IfFileExists{hyperref.sty}{\hypertrue}{\hyperfalse}
\hyperfalse  
\ifhyper\usepackage[naturalnames]{hyperref}\fi

\newif\ifdraft
\drafttrue
\numberwithin{equation}{section}
\numberwithin{figure}{section}

\newtheorem{theorem}{Theorem}
\numberwithin{theorem}{section}
\newtheorem{corollary}[theorem]{Corollary}
\newtheorem{lemma}[theorem]{Lemma}
\newtheorem{prop}[theorem]{Proposition}
\newtheorem{proposition}[theorem]{Proposition}

\def\eref#1{(\ref{#1})}
\newcommand{\R}{{\mathbb R}}
\newcommand{\C}{{\mathbb C}}
\newcommand{\U}{{\mathbb U}}
\newcommand{\N}{{\mathbb N}}

\def \E {{\bf E}}
\def \A {{\cal A}}
\def\ev#1{\mathcal{#1}}
\def \L {{\cal L}}

\def\p{\partial}
\def\gp{{\gamma'}}
\def\tgp{{{\tilde\gamma}'{}}}
\def\bgp{{{\bar\gamma}'{}}}

\def \eps {\varepsilon}

\def \P {{\bf P}}
\def\Bb#1#2{{\def\md{\bigm| }#1\bigl[\,#2\,\bigr]}}
\def\BB#1#2{{\def\md{\Bigm| }#1\Bigl[\,#2\,\Bigr]}}
\def\Bs#1#2{{\def\md{\mid}#1[\,#2\,]}}
\def\Pb{\Bb\P}
\def\Eb{\Bb\E}
\def\PB{\BB\P}
\def\EB{\BB\E}

\def\Es{\Bs\E}
\def\Pw#1/#2/{\widetilde\P_#1^#2}
\def\coupl{\mu}
\def \K {{\cal K}}

\def\ZZ{\hat Z}
\def \Disk {\U}
\def \bg {{\bar \gamma}}

\def \Prob {{\bf P}}

\def\closure#1{{\overline{#1}}}

\def\proofof#1{{ \medbreak \noindent {\bf Proof of #1.} }}
\def\proof{{ \medbreak \noindent {\bf Proof.} }}
\def\st{\, : \,}
\def \rv {\Xi}

\def \r {{\rho}}

\def\Im{\mathrm{Im}}
\def\Re{\mathrm{Re}}

\begin{document}

\maketitle

\begin {abstract}
We show that the intersection exponents for planar Brownian 
motions are analytic. More precisely, let $B$ and $B'$ be 
independent planar Brownian motions started from 
distinct points,
and define the exponent $\xi (1, \lambda)$ by
$$
\EB{\Pb{ B[0,t] \cap B'[0,t] = \emptyset
  \md B[0,t]}^{\lambda}}
\approx 
t^{-\xi(1, \lambda)/2}, \qquad t \to \infty.
$$
Then the mapping $\lambda \mapsto \xi (1, \lambda)$ is
real analytic in $(0,\infty)$.
The same result is proved for the exponents $\xi (k, \lambda)$
where $k$ is a positive integer.
In combination with the determination of
$\xi ( k, \lambda)$ for integer $k \ge 1$ 
and real  $\lambda \ge 1$  
in our previous papers, this gives the 
value of $\xi (k, \lambda)$ also for 
$\lambda \in (0,1)$ and the disconnection
exponents $\lim_{\lambda \searrow 0} \xi (k, \lambda)$.
In particular, it shows that
$\lim_{\lambda \searrow 0} \xi(2, \lambda ) = 2/3$ and concludes
the proof of the following result that had
been conjectured by Mandelbrot:
the Hausdorff dimension of the  outer boundary of $B[0,1]$
is $4/3$ almost surely.
\end {abstract}
 

\section {Introduction}

The goal of the present paper is to show that the intersection 
exponents for planar Brownian motions are analytic.
 Let $k \ge 1$ be a  positive integer 
and let $X^1, \ldots, X^k$ be independent Brownian motions in
the complex plane $\C$
started from $0$.
Let $Y, Y^1, Y^2, \ldots $ denote other independent planar Brownian motions
started from 
$1$, and let $\rv_t$ be the random variable (measurable with respect
to $X^1,\ldots,X^k$),
\[  \rv_t = \Pb{ Y[0,t]\cap (X^1 [0,t] \cup \cdots \cup X^k [0,t]) 
= \emptyset  \md X^1 [0,t] \cup \cdots \cup X^k [0,t]} . \]
Note that
\[  \Prob
\bigl[ (X^1 [0,t] \cup \cdots \cup X^k [0,t]) 
\cap (Y^1 [0,t] \cup \cdots \cup Y^p [0,t] ) = \emptyset 
\bigr] = \E[\rv_t^p] . \]
The intersection exponent $\xi(k,\lambda)$ is defined for $\lambda> 0$ by
\begin{equation}  \label{may7.1}
   \E[\rv_t^\lambda] \approx 
t^{-\xi (k,\lambda)/2}, \qquad  t \to \infty,
\end{equation}
that is,
$$
\xi(k,\lambda) :=-2\lim_{t\to\infty} \frac{\log\E[\rv_t^\lambda]}{\log t}
.$$
The existence of such  exponents follows easily from a subadditivity 
argument.
For a more detailed account of the definition and properties of these
exponents, we refer the reader to our earlier papers
\cite {LW1,LSW1,LSW2}.
Let us mention, however, that they are related to 
other critical exponents arising in statistical physics,
including those 
predicted by theoretical physicists 
for planar critical percolation and self-avoiding walks
(see references in \cite {LSW1}).

In \cite {LSW2}, the value of 
$\xi (k,\lambda)$ was determined for a large collection of pairs
$(k, \lambda)$.
In particular, it was shown that 
\begin{equation}
\label {apr5}
\xi(k,\lambda)
=
\frac{(\sqrt{24 k+1}+\sqrt{24\lambda+1}-2)^2-4}{48}
\end{equation}
holds for $k=1 , \lambda \ge 10/3$ and 
for $k=2, \lambda \ge 2$.
In \cite {LSW2s}, we then showed that 
(\ref{apr5}) holds
for all integer $k \ge 1$ and all real $\lambda\ge 1$.
The idea in the proofs is to compute the 
exponents associated to another conformally
invariant process (called stochastic Loewner evolution 
process and first introduced  
in \cite {S1}, see also \cite {S2}) 
and to identify them with the 
Brownian intersection exponents via a universality 
argument (introduced in \cite {LW2}).
This universality argument is not sufficient to 
derive the value of $\xi (k, \lambda)$ when $\lambda<1$. 
 
\medbreak

The
exponents $\xi (2, \lambda)$ 
are very closely related to 
the dimension and properties of the so-called outer boundary of a planar
Brownian path. 
The outer boundary, or frontier,
of the Brownian path $B[0,1]$ is the boundary of the unbounded
 connected
component of $\C\setminus B[0,1]$.  
The disconnection exponents $\eta_k$ 
are defined by
$$
\Pb{X^1 [0,t] \cup \cdots \cup X^k [0,t] \hbox { does not disconnect } 1 
\hbox { from } \infty}
\approx t^{-\eta_k/2}.
$$
It is easy to
show that $\eta_2 
= \xi (2, 0)$ if we use the
convention $0^0 = 0$ in the definition of $\xi (k, \lambda)$ when $ 
 \lambda= 0$. 
In \cite {Lfront} it was proved that
 the Hausdorff dimension of the  frontier of $B[0,1]$
is almost surely equal to $2 - \eta_2$.
Moreover, \cite {Lmulti} the multifractal spectrum of the frontier 
with
respect to harmonic measure is 
given in terms of the Legendre transform 
of the function $\xi (2, \lambda)$. 
In \cite {Lstrict}, it is also shown that 
$$
\lim_{\lambda \searrow 0}
\xi(k, \lambda) = \xi(k,0) = \eta_k\,.
$$

\medbreak

The main result of the present paper is the following:

\begin {theorem}
\label {main}
For all integers $k \ge 1$, the function $\lambda  \to \xi (k, \lambda)$ 
is real analytic in $(0,\infty)$.
\end {theorem}

This has the following consequences:

\begin {corollary}
Formula \eref{apr5} is valid for all positive integers $k$
and all non-negative real  $\lambda$.
In particular,
\begin {equation}
\label {disco}
\eta_k  = \frac { \left( \sqrt { 24 k + 1 } - 1 \right)^2 - 4 }{48}.
\end {equation}
\end {corollary}

\proof For $\lambda>0$, this follows
from \eref{apr5} and Theorem~\ref{main} by analytic continuation.
For $\lambda=0$, the result follows by
the continuity at $0$ proved in \cite{Lstrict}.
\qed

\begin {corollary}\label{BMdim}
The Hausdorff dimension of the outer boundary of a planar
Brownian path $B[0,1]$ is almost surely $4/3$.
\end {corollary}
\proof The case $k=2$ in~\eref{disco} gives $\eta_2=2/3$.
The corollary follows from this and the result from~\cite{Lfront}
saying that the dimension of the frontier is $2-\eta_2$.
\qed
\medbreak

Note that for this corollary, one does not need \cite{LSW2s},
since \eref{apr5} appears in \cite{LSW2} for $k=2$
and  $\lambda \ge 2$.

Corollary~\ref{BMdim} has been conjectured by Mandelbrot \cite {M},
based on simulations and the analogous
conjecture for self-avoiding walks.
 Nonrigorous arguments from 
theoretical physics 
involving quantum
gravity \cite {Dqg} also lead to this conjecture.
Before our series of papers
\cite {LSW1,LSW2,LSW2s}, it had been proved
\cite {Bal, BL2, Wecp,Lfront}
that the Hausdorff dimension of the outer boundary 
of a planar Brownian path is in the
interval $(1.01,1.48)$.
    
The Hausdorff dimension of other exceptional 
subsets of the planar Brownian curve can be described
in terms of disconnection exponents.
A point $z$ is a pioneer point of $B[0,1]$
 if there is some
time $t\in[0,1]$ such that $z=B_t$ is in the outer
boundary of $B[0,t]$.  It is shown in \cite{Lbuda} 
that the Hausdorff dimension of the set of pioneer points
is $2-\eta_1$ almost surely.  Consequently,~\eref{disco} gives

\begin{corollary}
The Hausdorff dimension of the set of pioneer points of a planar
Brownian path is almost surely $7/4$.
\qed
\end {corollary}

In the same way, one 
gets that the Hausdorff dimension of the 
set of double points of $B[0,1]$ that are 
also on the outer boundary  of $B[0,1]$
is $2 -  \eta_4
= ( 1 + \sqrt {97}) / 24 $ (which is not a rational 
number).
\medskip

To prove Theorem \ref{main}, we
show that for every $\lambda > 0$, the function
$x\mapsto \xi(k,x)$ can be extended to an analytic function
in a neighborhood of $\lambda$ in the complex plane. 
For notational ease,
we will restrict the proof to the case $k=2$; the proof for
other values of $k$ is essentially the same.
Our proof has similarities with the proof that the free energy of
a one-dimensional Ising model with exponentially
decaying interactions is an analytic function (see \cite{R}).
We shall show that $e^{-\xi(2,x)}$ is the leading eigenvalue
of an operator $T_x$
on a space of functions on pairs of paths.
A special norm will be chosen such that on the space of functions
with finite norm, the dependence of $T_x$ on $x$ is analytic,
and the leading eigenvalue is an isolated
simple eigenvalue.
 It then follows by
a standard result from operator theory
(see, e.g., \cite {DS})
that $\xi(2,x)$ is an analytic function of $x$.

Let us outline the argument in the paper.  The Banach
spaces and the operators are defined in Section~\ref{s.setup}.
The operators act on spaces of functions of
pairs of paths from the origin
to the unit circle.
The function spaces consist of functions
with the property that their dependence  
on the behavior of the paths near the origin decays
exponentially.
These 
spaces are reminiscent of the spaces defined in \cite{R}, but the
precise definition in this paper is new.  In Section
\ref{prevsec} we review facts about $\xi(2,\lambda)$ from
\cite{Lstrict} which are needed in the proof.
Analyticity
of the operator is proved in Section~\ref{analsec}, using a 
coupling argument.
The existence of the
spectral gap is established
in Section  \ref{specsec}.  The main tool to
show that the 
eigenvalue is isolated is also a coupling
--- but this time a coupling of weighted Brownian paths. 
  A similar coupling
was used in \cite{BFG} for a one-dimesional
Ising-type model,
and 
 such a coupling was
 first used  in \cite{Lhalf} for weighted
Brownian paths.

We end the introduction with
 a few brief words regarding notation.
In this paper,  $c,c',u $, etc., 
 denote constants, whose values may  change from line to line,
 while 
$C$ and  $c_0, c_0', v_0$, etc.\ will denote constants 
whose value will not
change.  The values of these constants
will depend on $\lambda$.
 The notation
$f(t)\approx g(t)$ means $\log f(t)/\log g(t)\to 1$ 
as $t\to\infty$, while $f(t)\asymp g(t)$ means
that there is a constant $c>0$ such that $c^{-1}\le f/g\le c$.
The unit disk $\bigl\{|z|<1\bigr\}$ in $\C$ is denoted by $\U$.

\section {The operator}\label{s.setup}

Let $\Gamma_0$ denote the set of all 
continuous paths
$\alpha:[0,1]\to \closure{\U}$ such
that $\alpha(0)=0$, $|\alpha(1)|=1$
and $0<|\alpha(t)|<1$ for $t\in (0,1)$.
We identify  two paths if one can
be obtained from the  other by
increasing reparameterization,
and endow $\Gamma_0$ with the  metric
$$
d(\alpha,\beta):= \inf_{\phi}\sup_{t\in[0,1]}
\left|\alpha(t)-\beta\bigl(\phi(t)\bigr)\right|\,,
$$
where $\phi$ ranges over all increasing
homeomorphisms $\phi:[0,1]\to[0,1]$.
Let $\Gamma\subset\Gamma_0\times\Gamma_0$
denote the set of $\gamma
= (\alpha, \beta)$
such that there exists a unique connected
component $O=O(\gamma)$ of $\U \setminus ( \alpha \cup \beta)$
with $0 \in \partial O$ and $\partial \U \cap \partial O 
\ne \emptyset$.
 Let $\A$ denote  the set of bounded
measurable  functions $f:\Gamma\to\C$
and $\| f \| = \sup_{\gamma \in \Gamma} | f(\gamma)|$. 
Here, measurability is with respect to the Borel
sets from the metric on $\Gamma$ (and this is the sole use of
this metric in this paper).

We are interested in  functions $f \in \A$ that 
depend little  on the part of $\gamma$
near  the origin.
For $\alpha\in\Gamma_0$ and $m\ge 0$, let
$\alpha_m$ be the arc of $\alpha$ after its first
point in the circle $e^{-m}\p\U$ of radius $e^{-m}$ around $0$,
and for $\gamma=(\alpha,\beta)\in\Gamma$
set $\gamma_m=(\alpha_m,\beta_m)$.
If $k < j$,
we say that the path $\alpha \in \Gamma_0$  
has no downcrossing from $e^{-k}$
to $e^{-j}$
if $\alpha_k \cap e^{-j} \U = \emptyset$, in other words,
if $\alpha$ does not touch $e^{-j} \U$ after its first visit to
$e^{-k} \p \U$.
We say that $\gamma = ( \alpha, \beta)$ 
has no downcrossing if both $\alpha$ and $\beta$ have no
downcrossing.
Let $\ev Y_m$ be the set of all $\gamma \in \Gamma$
such that 
 for all $k \in  [0, 11m/12]$,
 $\gamma$ contains no downcrossing from $e^{-k}$
to $e^{- k - m/12}$.
Let  $\ev X_m$ be the set
of $( \gamma, \gamma') \in \Gamma^2$ such that 
\begin {itemize}
\item $( \gamma, \gamma') \in \ev Y_m \times \ev Y_m$,
\item $\gamma_m = \gamma_m'$, and
\item $O(\gamma) \cap \partial \U
= O ( \gamma') \cap \partial \U$.
\end {itemize} 

For all $f \in \A$ and $u>0$, let
$$
\| f \|_u
:= 
\max \Bigl\{
\| f \|
,\,
\sup\bigl\{
e^{mu} | f(\gamma) - f (\gamma') |
: m=1,2,\dots,\, (\gamma,\gamma')
\in \ev X_m
\bigr\} \Bigr\}.
$$
This is a norm on the Banach space 
$\A_u := \{ f \in \A \, : \, \|f\|_u< \infty
\}.$
Let ${\cal L}_u$ denote  the Banach space of
continuous linear operators from $\A_u$
to $\A_u$ with the operator norm
$$
N_u (T) := \sup_{\| f \|_u = 1} \| T(f) \|_u\,.
$$

For every $\gamma \in \Gamma$, let
$\hat \alpha$
and $\hat \beta$ be
 independent planar Brownian motions 
 started, respectively, from the endpoints of
 $\alpha$ and $\beta$ 
 on the unit circle. Let
$\hat \alpha^n$  denote the path 
$\hat \alpha$  stopped when it 
hits $e^n\p\U$;
let  $\tilde \alpha^n$ be the path from $0$ to 
$e^n\p\U$ obtained by
concatenating $\alpha$ and $\hat \alpha^n$ ;
and let $\bar\alpha^n:=e^{-n}\tilde\alpha^n$.
Define $\hat \beta^n,\tilde \beta^n,$
and $\bar \beta^n$ similarly, and let
$\hat\gamma^n:=(\hat\alpha^n,\hat\beta^n)$,
$\tilde \gamma^n := ( \tilde \alpha^n, \tilde \beta^n)$,
$\bar \gamma^n := (\bar\alpha^n,\bar\beta^n)$.
We will often omit the superscript $n$ when $n=1$.

Define the event ${\cal E}_n := \{ \bar \gamma^n  \in \Gamma \}$.
Note that almost surely this event is satisfied as long as 
 $\tilde \gamma^n$ does not disconnect $0$ from $\infty$
(since  $0 \notin \hat \gamma^n$ a.s.).

For every $\gamma \in \Gamma$, consider the $h$-process
$B$ started at $0$ and conditioned to first leave $O= O (\gamma)$
in $\partial \U \cap \partial O$.  
Let us say a few words about how this
process is defined: it can 
be viewed as
 the
limit as $z \in O$ tends to $0$ of Brownian motion
starting from $z$ conditioned to leave
$O$
in $\partial \U \cap \partial O$.  It is well-defined
since $0$ is a simple boundary point of $O$ for
$\gamma \in \Gamma$.
For instance, if we map  
 conformally 
the strip  $(0,1)
\times \R $
onto $O(\gamma)$
 taking  $\{1\}  \times \R$ 
onto $\partial \U \cap \partial O$
 and the origin to the 
origin, then the $h$-process in $O(\gamma)$ is obtained (after
time-change, but we will actually only 
use the paths of the $h$-processes) 
as the image under the conformal map of  
the process $X + i Y$ in the strip, where $X$ is a 
three-dimensional Bessel process
started from 0  and stopped when 
it hits 1 (i.e., the limit 
when $\eps \to 0$
of one-dimensional 
Brownian motion started from $\eps>0$ conditioned to hit 1
before 0),
and $Y$ is an independent one-dimensional Brownian motion
started from 0 (stopped at the same time).

Attach to the endpoint of $B$ on $\partial \U$
an independent Brownian motion $\hat B$, 
and define the paths $\hat B^n$ and  $\tilde B^n$ as before.
The path $\tilde B^n$ consists of two parts: The $h$-process
(up to its hitting time of $\partial \U$), and the 
(non-conditioned) Brownian motion $\hat B^n$.
Let
\[
Z_n  = Z_n (\tilde \gamma^n) :=
\Pb {  \tilde B^n
  \cap \tilde \gamma^n  = \emptyset \md \tilde\gamma^n},
\]
and $Z=Z_1$.  
Note  that  although
$B \cap \gamma  = \emptyset$ almost surely,
 it can happen that 
$B \cap \hat \gamma^n \ne \emptyset$ with positive
probability.  That is,
$ \tilde B^n  \cap
\tilde \gamma^n= \emptyset $
can fail in two ways: it may happen
that $\hat B^n\cap\tilde\gamma^n\ne\emptyset$,
and it may also happen that $B\cap \hat\gamma^n\ne\emptyset$.
Note also that $Z_n \ne 0$ if and only if ${\cal E}_n$ 
occurs.
We define
$\psi_n = - \log Z_n$ and $\psi=
- \log Z$ (with $- \log 0 = \infty$).

For $n,\lambda>0$ define the linear operator $T_\lambda: \A \rightarrow
\A$ by
\[ T_\lambda^n f(\gamma) :=
\EB{ f(\bar\gamma^n) \exp\bigl(-\lambda \psi_n \bigr) }
=
\EB{ f(\bar\gamma^n) Z_n^\lambda }
 . \]
and let $T_\lambda = T_\lambda^1$.  The expectation is
over the randomness in $\tilde \gamma^n$.  Note that 
$T^{n+m}_\lambda = T^n_\lambda T^m_\lambda$,
so the notation is appropriate.
Also, there is no need to restrict to real $\lambda$; this 
defines $T_z^n$ for complex $z$ with $\Re(z) > 0$.
We will prove the following:

\begin {prop}
\label {OK}
\begin {itemize}
\item [{(i)}]
For all real $\lambda>0$, there exist an $\eps >0$ and a $v = v(\lambda) >0$
such that for all $ u\in(0,v)$, $z \mapsto T_z$
is an analytic function from 
 the disk $\bigl\{z \st |z-\lambda|<\eps\bigr\}$ into $\L_u$. 
\item [{(ii)}]
For all  real  $\lambda>0$, there exist an $\eps'>0$, 
and a $u\in(0, v (\lambda))$
such that the spectrum of $T_\lambda$
in $\L_u$ is the union of 
the simple eigenvalue 
$e^{-\xi(2, \lambda)}$ and  a subset of 
the disk $ (1-\eps')e^{-\xi(2, \lambda)}\U$
\end {itemize}
\end {prop}

\proofof{Theorem \ref{main} {\rm (assuming  Proposition \ref{OK})}}
The proposition implies that $e^{-\xi(2,\lambda)}$ is an isolated
simple eigenvalue of $T_\lambda$, for all $\lambda>0$.
By Theorem VII.6.9 in \cite {DS}, it follows 
that for all $\lambda>0$, there exists $\eps>0$ such that 
$x \mapsto \xi (2, x)$ can be extended analytically to the disk
$\bigl\{z \st |z-\lambda|<\eps\bigr\}$.
Hence, there exists a neighborhood ${\cal N}$ of the half-line $(0, \infty)$
such that $z \mapsto \xi (2, z)$ is  analytic in ${\cal N}$,
proving the theorem in the case $k=2$.
The proof for other $k$ is the same. \qed

The proofs of Proposition~\ref {OK}(i) and (ii) both rely 
on  coupling arguments.  The proof of (i) in Section \ref {analsec}
 uses a coupling 
of the $h$-processes $B$ and $B'$ associated to 
two pairs of paths $\gamma$ and $\gamma'$  
(when $(\gamma,\gamma') \in \ev X_m$).
In the proof of (ii) (Section \ref {specsec}),  
we couple the extensions
$\hat \gamma^n$ and $\hat {\gamma'}^n$ 
associated to $\gamma$ and $\gamma'$ defined under a 
weighted probability measure.

\section {Previous results on $\xi(2,\lambda)$}  \label{prevsec}

In this section, we very 
briefly review some important facts 
about the intersection exponent $\xi(2,\lambda)$
that will be useful.
The results here were derived in \cite{Lstrict} and
apply to dimensions $2$ and $3$.  Since some of
the arguments are simpler when one considers only the
plane,  we plan to release  
detailed self-contained proofs of all these facts 
((\ref {up2}), (\ref {sep}) and (\ref {tildec}))
very shortly \cite {LSWup2}.

\subsection{Estimates up to constants}

Let $x,y\in\p\U$, and
let $W=(X,Y)$ denote a pair of planar Brownian motions started
at $x,y$, respectively.
Let $B$ denote another independent
 planar Brownian motion 
started from $b\in\p\U$.
Let $B^n$ be the path $B$ until it reaches
 the circle $e^n\p\U$ for the first time, and similarly
define $X^n$ and $Y^n$.
Let $W^n =  X^n \cup Y^n$,  
$$
{\ZZ}_n = {\ZZ}_n (W^n) 
:= \sup_{b\in\p\U} \Pb { B^n \cap W^n = \emptyset \mid W^n },
$$
and 
$$
q_n=q_n(\lambda):=\sup_{x,y\in\p\U}
\Eb{\ZZ_n^\lambda}.
$$
Then $q_{n+m}\le q_nq_m$ holds, by the strong Markov property
for $X$, $Y$ and $B$.  Hence,
the limit
$
\xi := -\lim_n \log(q_n) / n = -\inf_n   \log(q_n) /n,
$
exists by subadditivity and 
\begin {equation}
\label {lower} 
q_n \ge e^{-\xi n}\,.
\end {equation}
It is not difficult to verify that, in fact, $\xi=\xi(2,\lambda)$,
as defined in the introduction.

An opposite estimate also holds; that is,
there exists a  constant $c_1= c_1(\lambda)>0$, such that 
\begin {equation}
\label {up2}
q_n  \le c_1 e^{- \xi n}.
\end {equation}
This result is a variant  of Theorem 2.1 in \cite {Lstrict}.
The fact that $q_n \asymp e^{-\xi n}$ (which is 
a more precise statement  than 
the definition of $\xi$, $q_n \approx e^{-\xi n}$)  
has been instrumental in showing that $2- \xi$
correspond to Hausdorff dimensions 
of various subsets of the planar Brownian
curve \cite {Lmulti}.

\subsection{Separation and the functions $R_n$}
\label{Ssep}

 Let  $\Gamma^+$ be the set of $\gamma = (\alpha,\beta)
\in \Gamma$ such that
\[  (\alpha_{1/2},\beta_{1/2}) \subset 
\bigl\{  e^{\r+i \theta} \ : \ \r \in [-1,0] ,
\, \theta \in (\pi /2, 3\pi/2)\bigr\}, \]
and
$$ \bigl\{ e^{\r+ i \theta} \ : \ \r \in (-1/2, 0) ,\, 
\theta \in (-\pi/2, 
\pi/2)
\bigr\} \subset O (\gamma)
.$$
The Separation Lemma
\cite[Lemma 4.2]{Lstrict}
states that there exists a $c_2$ such that for all $\gamma \in \Gamma$,
\begin{equation}  \label{sep}
  \Eb{\exp(-\lambda \psi)\, 1_{\{\bar \gamma \in \Gamma^+\}} }
   \geq c_2 \,\Eb{\exp(-\lambda \psi)}. 
\end{equation}
It 
is a kind of boundary Harnack principle
for the operator $T$. 
This type of result was important in the derivation 
of (\ref {up2}).

Define for all $n \ge 1$, $\lambda>0$ and $\gamma\in\Gamma$,
$$
R_n(\gamma)  = R_{n,\lambda}(\gamma) :=
 e^{n \xi} \E [ e^{- \lambda \psi_n} ]
= e^{n\xi} (T^n_\lambda 1)(\gamma).$$ 
By (\ref{up2}), $R_n$ is uniformly bounded in $n$ and 
$\gamma$. Note also that, by the strong Markov property,
and~\eref{up2}
\begin {equation}
\label {n<1}
 R_n (\gamma) = e^{n\xi} (T^n_\lambda 1)(\gamma)
     \leq e^{-n \xi} q_{n-1} T_\lambda 1(\gamma)
      \leq c R_1(\gamma) . 
\end {equation}
On the other hand, 
one can prove that 
\begin {equation}
\label {tildec}
\tilde c:= 
\inf_{n\ge 0}\inf_{\gamma\in\Gamma^+}
R_n(\gamma) >0.
\end {equation}
This is, loosely speaking, due to the fact that a positive fraction 
of the extensions $\hat B^n$ and $\hat \gamma^n$ leave quickly 
the neighborhood of $\U$, so that they do not really feel
the influence of $\gamma$ and $B$. As $\gamma \in \Gamma^+$, 
the starting points of $\gamma^n$  
and a positive fraction of the starting 
points of $\hat B^n$ 
 are not so far from being optimal,
so that $q_n$ is within a constant multiple of
 $\Es{  e^{-\lambda \psi_n} }$.

Using the strong Markov property
and~\eref{sep}, then
applying (\ref {tildec}) to 
$R_{n-1} ( \bar \gamma)$, gives
for all $\gamma \in \Gamma$ and $n \ge 1$,
\begin {equation*}
  R_n(\gamma) \geq e^{n\xi} \Eb{e^{- \lambda \psi_n}\,
      1_{\{\bar \gamma \in \Gamma^+\}}}
      \geq e^{n \xi}\, \tilde{c}\, e^{-(n-1)\xi}
             \Eb{e^{-\lambda \psi} 
\,1_{\{\bar \gamma \in \Gamma^+\}}
  } \geq 
  \tilde{c}\, c_2 \,R_1(\gamma).
\end {equation*}
Combining this with (\ref {n<1}) shows
that  there is a $c_0$ such for all $n,n' \geq 1$
and  $\gamma \in
\Gamma$,
\begin{equation}
\label{use}
R_n(\gamma) \leq c_0 R_{n^\prime} (\gamma) . 
\end{equation}
In \cite{Lstrict} it was shown that the limit
$ R(\gamma) = \lim_{n \rightarrow \infty} R_n(\gamma)$ 
exists.
We will rederive this result in this paper and simultaneously
improve the rate of convergence to the limit.

\section {Analyticity}  \label{analsec}

The goal of this section is to
 prove  Proposition \ref {OK}(i).

\subsection{Coupling $B$ and $B'$}
\label {S.B}

Let $\gamma,\gamma'\in\Gamma$.
Let $B$ be the $h$-process associated with $\gamma$,
and let $B'$ be the $h$-process associated with
$\gamma'$.  In this section, we show that there is
fast decay for the dependence of $B$ on $\gamma$.
More precisely, we prove the following proposition.

\begin{prop}
\label {p.decay}
There exist a $ c> 0$ and a $v_0>0$ such that  for all $m\ge 1$,
if $(\gamma,\gamma')\in\ev X_m$, then
one can couple $B$ and $B'$ in such a way that 
$$
\Pb { B \setminus e^{-m/2} \U 
\ne B' \setminus e^{-m/2} \U } \le c\: e^{-v_0 m }.
$$
\end{prop}

By coupling $B$ and $B'$, we mean that it is possible to define 
$B$ and $B'$ on the same probability space 
in such a way that the law of $B$ (resp. $B'$) is that 
of the $h$-process associated with $\gamma$ (resp. $\gamma'$).
A reference for facts about the coupling method is \cite {Li}.
 
This result
 actually holds for all $(\gamma,\gamma^\prime)
\in \Gamma$
with $\gamma_m= \gamma^\prime_m$ and $O(\gamma) \setminus e^{-m}\Disk
= O(\gamma^\prime) \setminus e^{-m} \Disk$, but the proof is
more complicated.  Proposition \ref{p.decay}
 will be sufficient for our purposes.

It is easy to verify that the processes $B$ and $B'$ satisfy the strong 
Markov property.  This will be very useful in the following.

In preparation for the proof of Proposition \ref {p.decay}, we 
first focus on Brownian motions in  the half-infinite rectangle
\[   J = \{x + iy: 0 < x < \pi , 0 < y < \infty \} , \]
and then apply these
results to conditioned Brownian motions 
in $O$. There is a conformal transformation taking $O(\gamma)$ to
$J$ which takes the origin to infinity and $O(\gamma)
\cap \partial \Disk$ to $[0,\pi ]$, so that  
conditioned Brownian motions 
in $J$ can be conformally transformed (up to a time-change) to
conditioned Brownian motions in $O(\gamma)$.

\begin {lemma}
 \label{lemma.apr26}
(i) 
There exists a constant $c$ such that for all 
$z $ in $J$ such that $ \Im (z) \ge 1 $ and  for all $y_0 \ge 1$, 
if $\tilde B$ denotes a 
Brownian motion started from $z$, and conditioned to 
leave $J$  on the set $[0,\pi]$, then
\begin {equation}
\label {firstfact}
\Pb{   
\tilde B \hbox { hits } \{w \st \Im (w) \ge \Im (z) + y_0  \}
\hbox { before it leaves } J   }
 \leq  c e^{-2y_0} . 
\end {equation}
\noindent(ii)
There exists a constant $c$ such that for 
 every $\tilde m\ge 0$, every $\tilde n\ge 1$ and every $z,z' \in J$
with $\tilde m+\tilde n \leq \Im(z) \le \Im(z^\prime)$,
one can find a coupling of $\tilde B$ and $\tilde B^\prime$,
Brownian motions
conditioned to leave $J$ at $[0,\pi ]$,
  starting at $z$ and $z^\prime$
respectively,  such that
\begin {equation}
\label {cou}
\PB{
    \tilde B \cap \{ w \st \Im(w) \leq 
 \tilde m  \} \ne
   \tilde B' \cap \{ w \st \Im(w) \leq \tilde m  \} }
\le
c e^{-\tilde n} . 
\end {equation}
\end {lemma}

\proofof{Lemma \ref{lemma.apr26}}
Let  
$\tilde \tau$ denote the first time a Brownian motion $\tilde
X$ (started from $\tilde X=z$ 
under the probability measure $\P^z$)
leaves $J$, and for all $y_0>0$, let $\tilde \sigma_{y_0}$
denote the first time at which $\Im ( \tilde X) $ hits $y_0$.

Let $h_z(s),\, z \in J ,\, s \in (0,\pi)$ be the density of
$\tilde X_{\tilde \tau} 1_{\{\tilde \tau = \tilde \sigma_{0}\}}$,
where $\tilde X_0=z=x+iy$.
  This density can easily be found,
for instance, by separation of variables;  
in particular, it is easy to show that 
 the density satisfies for all $y \geq 1$
\begin {equation}
\label {h=}
h_z (s) = \frac{2}{\pi}\sin(x) \sin(s)  e^{-y}
      \,\bigl ( 1+ O(e^{-y})\bigr)    ,
\end {equation}
where the constant implicit in the $O(e^{-y})$ notation
does not depend on $x,s$, or $y$.
By reflecting $z$ in the line $\Im(w)=y+y_0$, it follows that
\begin {eqnarray*}
 \Prob^z 
[ \tilde \sigma_{y+y_0} < \tilde \tau = \tilde \sigma_0 ]
&\le & \Prob^{z+2iy_0} 
[ \tilde \sigma_0 = \tilde \tau ]
\\
& = &
 \frac {4}{\pi} \sin (x) e^{-y-2 y_0} \bigl ( 1+ O (e^{-y})\bigr)
\\
& = & 
 \bigl(1 + O (e^{-y})\bigr )\, e^{-2 y_0} \,
\Prob^z [ \tilde \tau = \tilde \sigma_0],
\end {eqnarray*}
and (\ref {firstfact}) follows.
\medskip

Let $\tilde z=i\tilde m+\tilde s$ be the first point
on the segment $i\tilde m+[0,\pi]$ which $\tilde B$ hits,
and let $\tilde z'=i\tilde m + \tilde s'$ be the
corresponding point of $\tilde B'$.
A straightforward consequence of the 
strong Markov property shows that 
the density of $\tilde s$ is 
$$
\tilde h_{z} (s) :=
\frac{
h_{z-i\tilde m}(s)\int_{0}^{\pi}h_{s+i\tilde m}(t)\,dt}
{\int_{0}^{\pi} h_z(t)\,dt}\,,
$$
and a similar expression holds for $\tilde h_{z'}$,
the density of $\tilde s'$.
Using~\eref{h=}, this gives
$$
\tilde h_z(s)/\tilde h_{z'}(s) = 1+O\bigl(e^{-\tilde n}\bigr).
$$
Therefore, we may couple the parts of $\tilde B$ and $\tilde B'$
until their first hit of $i\tilde m+[0,\pi]$
so that $\P[\tilde z\ne\tilde z']\le O (e^{-\tilde n})$.
On the event $\tilde z=\tilde z'$, we may continue
$\tilde B$ and $\tilde B'$ as the same path.
This gives \eref {cou}.
\qed

\begin{figure}
\SetLabels
(.63*.8)$O(\gamma)$\\
\R(.2*.5)$\alpha$\\
\L(.83*.4)$\beta$\\
\B(.45*.72)$A(\gamma,k)$\\
\R\B(.46*.42)$A(\gamma,j)$\\
\B\L(.5*.25)$0$\\
\L(.7*.05)$e^{-k}\p\U$\\
\L(.59*.15)$e^{-j}\p\U$\\
\L(.95*.4)$\p\U$\\
\endSetLabels 
\centerline{\AffixLabels{\includegraphics*[height=3.5in]{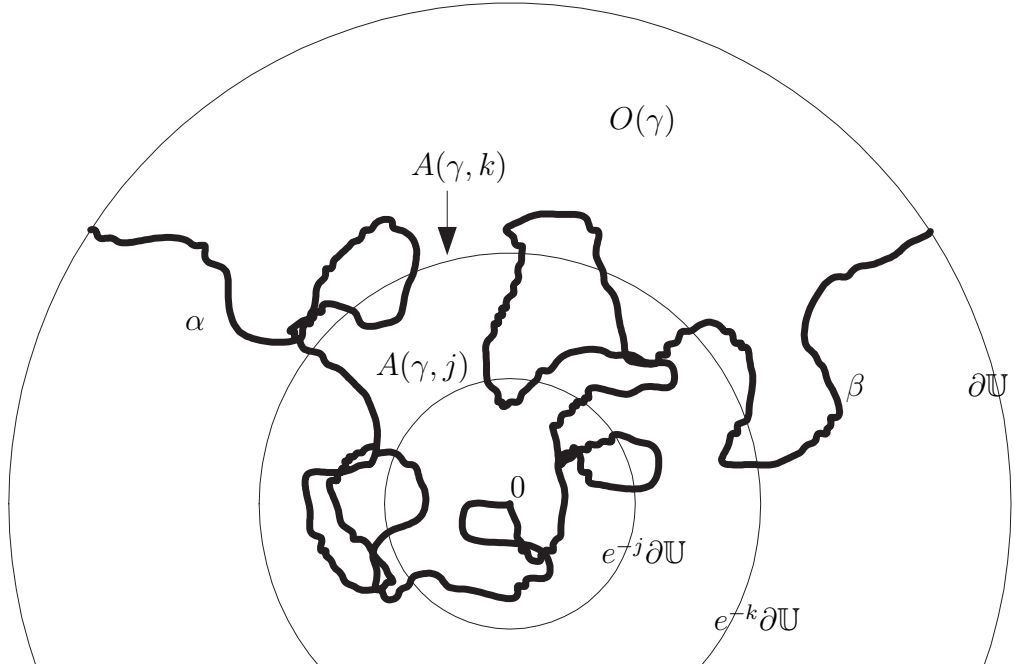}}}
\caption{\label{top}The separating arcs.}
\end{figure}

\proofof{Proposition \ref{p.decay}}
Prior to the core of the proof, there is a need
for some preliminary topological preparations. 
Suppose that $\gamma \in \Gamma$ and $k>0$.
Then, there exists at least
one open subarc 
$A = A(\gamma,k)$ of $ e^{-k}\p \U \cap O(\gamma)$ with
the property that the removal of $A$ from $O(\gamma)$
 disconnects zero
from $\partial \Disk$ in $O(\gamma)$.
See Figure~\ref{top}.
  There may be many such arcs; if so,
choose $A=A(\gamma ,k)$ 
as being the one closest to the origin, in the sense that every path
from the origin to $\partial \Disk$ in $O(\gamma)$
  goes through $A$ before any
other such disconnecting arc. 
The arc $A$ divides $O(\gamma)$ into
two components; let $O_+ (\gamma, k)$ denote the
 component whose
boundary 
intersects $\partial \Disk$, and let $O_- ( \gamma, k )$ denote the 
component whose boundary contains $0$.

As in Section \ref{s.setup}, let $\gamma_n = (\alpha_n,\beta_n)$
be the part of $\gamma$ starting
from the first visits to the circle
of radius $e^{-n}$.
 Note that $\gamma_k$ contains paths connecting the endpoints of
$A(\gamma,k)$ to $\p\U$.
Suppose for a moment that $j>k$ and
$\gamma$ has no downcrossing from $e^{-k}$ to $e^{-j}$.
Let $\tilde O_k(\gamma)$ denote the connected component of 
$\Disk \setminus \bigl( \gamma_k \cup A ( \gamma , k )\bigr ) $
that contains $O_+(\gamma,k)$.
Since $\gamma$ has no downcrossing from $e^{-k}$ to
$e^{-j}$, we have $\p \tilde O_k(\gamma)\cap e^{-j}\U=\emptyset$,
and therefore $\tilde O_k(\gamma)\cap e^{-j}\U=\emptyset$.
In particular, 
$$
O_+(\gamma,k)\cap e^{-j}\U=\emptyset.
$$
An entirely similar argument shows that
$$
O_-(\gamma,j)\subset e^{-k}\U
$$
if $\gamma$ has no downcrossing from $e^{-k}$ to $e^{-j}$.
In this case, one has to consider the connected component
that contains $O_-(\gamma,j)$
of the complement of the union
of $A(\gamma,j)$ with the parts of $\gamma$ until the last
visit to the circle $ e^{-j}\p\U$.
    
Now take $(\gamma,\gamma') \in \ev X_m$, as in the proposition.
Let $k \in[0, 11m/12]$;
so that $\gamma$ and $\gamma'$ have 
no downcrossing from $e^{-k}$ to $e^{-k-m/12}$.
Observe that $\tilde O(\gamma,k)=\tilde O(\gamma',k)$,
since both can be characterized as the largest
domain which does not contains zero,
has $\p O(\gamma)\cap\p\U=\p O(\gamma')\cap\p\U$
on its boundary, and is
bounded by $\gamma_k=\gamma'_k$ together with
an arc of $e^{-k}\p\U$.  This gives, 
$A(\gamma , e^{-k}) = A (\gamma', e^{-k})$.
Set
\begin {eqnarray*}
 A_1 &:=&
 A(\gamma,e^{-11m/12}) 
,\\
 A_2 &:=&  A(\gamma,e^{-10m/12}),
\\
A_3 &:=&  A(\gamma,e^{-9m/12}) ,
\\
A_\infty &:=&
 \partial O(\gamma) \cap \partial \Disk ,
\\
O_+ &:=&  O_+ (\gamma , e^{-11m/12}) ,
\\
O_- &:=&  O_- (\gamma, e^{-9m/12}).
\end {eqnarray*}
Except for $O_-$,
these are the same as the corresponding objects for $\gamma'$.
The extremal distance from $A_1$ to $A_2$ in $O(\gamma)$
is bounded from below by 
the extremal distance in $\C$ from $e^{-11m/12}\p\U$ to
$e^{-10m/12}\p\U$, which is $m/(12\cdot 2\cdot \pi)$.
Let $\phi$ be the conformal map from $O(\gamma)$
onto $J$ which takes $0$ to $\infty$ and takes $A_\infty$
onto $[0,\pi]$.  Note that $\phi(A_j)$ is a path joining
the two lines $\{\Re(w)=0\}$ and $\{\Re(w)=\pi\}$.
By conformal invariance of extremal length,
it follows that there is a constant $v>0$ such that
$$
\forall z_1\in \phi(A_1)\,
\forall z_2\in \phi(A_2)\qquad
\Im(z_1)-\Im(z_2)\ge v m-1/v\,.
$$

It follows immediately from Lemma \ref{lemma.apr26}.(i),
the strong Markov property and conformal 
invariance that the probability that $B$ hits $A_1$ after
its first hit to $A_2$ is bounded by $c e^{-vm}$.
The same holds for $B'$.
Note that if these processes do not hit $A_1$ again, then 
they stay in $O_+$.

To construct the coupling of $B$ and $B'$, let
them evolve independently until they 
first hit $A_2$ (at some random points $z$ and $z'$, respectively).   
The law of $B$ (respectively $B'$) after that time
is that of Brownian motion in $O(\gamma)$ (resp., $O(\gamma')$)
conditioned to exit $O(\gamma)$ (resp., $O(\gamma')$)
 in $A_\infty$.
 Each of these laws are at distance at most 
$ce^{-vm}$ from 
the law of Brownian motion
starting from $z$ and $z'$,
respectively, in 
$O_+$ conditioned to exit $O_+$ in $A_\infty$,
where distance is in the sense of the measure norm; that is, $L_1$-distance.
Hence, it suffices to show that for all 
$z, z' \in A_2$, it is possible to 
couple two Brownian motions
started from $z$ and $z'$, respectively,  conditioned
to exit $O(\gamma)$ in $A_\infty$, in such a way that they agree 
after their first hitting of $A_3$ with probability at least $1-ce^{-v'm}$.
(Note that before their first hitting of $A_3$, they are in 
$O_-  \subset e^{-8m/12} \Disk$.)
By conformal invariance, one
can use the coupling described in Lemma \ref{lemma.apr26}
with 
$$
\tilde m:=\sup\Im\bigl(\phi(A_3)\bigr),\qquad
\tilde n:=\inf\Im\bigl(\phi(A_2)\bigr)-\tilde m\,.
$$
By the conformal invariance of extremal length, as explained above,
$\tilde n/m$ is bounded away from zero, when $m$ is large.
Small values of $m$ can be handled by adjusting the constant $c$
in the statement of the proposition.
This completes the proof of Proposition~\ref{p.decay}.
\qed

\subsection {Exponential decay}

If $\gamma\in\Gamma$ and
$\gamma^\prime \in \Gamma$  have the same terminal point (which
will be the case, in particular, if $\gamma_m = \gamma^\prime_m$
for some $m >0$), then we will attach the same Brownian motions
$(\hat{\alpha},\hat{\beta}) = \hat \gamma = \hat \gamma'$
 to $\gamma$ and $\gamma^\prime$ in defining
$\tilde \gamma^n$ and $\tgp^n$. 
In this case,
write $\psi_n$ as shorthand for
$\psi_n(\tilde\gamma^n)$ and $\psi'_n$
as a shorthand for $\psi_n(\tgp^n)$.
Similarly, $Z_n'=Z_n(\tgp^n)$.

\begin {prop}
\label {ok3}
For every $\lambda >0$,
there exist $c>0$ and $v_1 > 0$, such that if $n\ge 0$, 
$m \ge 1$, and
$(\gamma, \gamma') \in \ev X_m$, 
then ,
$$
\Eb{\bigl| e^{-\lambda \psi_n} -
e^{- \lambda\psi_n'} \bigr|}
\le c e^{- \xi n} e^{-v_1 m}.
$$
\end {prop}
 
The proof of this proposition relies heavily on 
Proposition \ref {p.decay}.
An easy estimate on the disconnection
exponent will also be needed.
Let $X$ denote a planar Brownian motion started from the unit circle.
As before, for $\r>0$, let $X^\r$ denote the part of $X$ until
its first hit of the circle $e^\r\p\U$.
Let $v_0'>0$ be such that the probability that
$X^1$ disconnects the origin from infinity is 
$1- e^{-v_0'}$. Let $c_0' := e^{v_0'}$.
The strong Markov property immediately implies
that for all $\r >0$, 
\begin {equation}
\label {discon}
\Pb{ X^\r
\hbox { does not disconnect the origin from 
infinity}}
\le e^{v_0' (1-\r)} = c_0' e^{-v_0'\r}.
\end {equation}

\medskip

\proofof{Proposition \ref {ok3}} 
Suppose that $m \ge 1$
and that $(\gamma ,\gamma') \in \ev X_m$.
Couple the $h$-processes $B$ and $B'$ as in Proposition 
\ref {p.decay}.
When the coupling succeeds, i.e.,
when the event
$ {\cal H}
:= \{ 
\:B \setminus e^{-m/2} \U 
= B' \setminus e^{-m/2} \U \: \} 
$
holds, then
we attach to $B$ and $B'$ the same Brownian path $\hat B$.

\medskip

\noindent
{\sc Step 1.}
Define the events: 
$$\ev U (\hat \alpha) :=
 \{ \hat \alpha^n \cap e^{-m/2} \U = \emptyset \}
, \qquad 
\ev U(\hat \beta) :=
 \{ \hat \beta^n \cap e^{-m/2} \U = \emptyset \}
$$
and $\ev U=\ev U(\hat \gamma) :=\ev U(\hat \alpha )
\cap \ev U (\hat \beta)$.
When $\ev U(\hat \alpha)$ is not
satisfied, then   
 $\hat \alpha$ hits the circle  $e^{-m/2}\p\U$
before the circle $e^n\p\U$.
Let $\sigma$ be the first time at which $\hat\alpha(t)\in e^{-m/
2}\p\U$,
and let $\tau$ be the first time $t>\sigma$ such that
$\hat\alpha(t)\in\p\U$.
Let $\eta^n$ be the path $\hat \alpha$ after $\tau$ and until its
first hit of the circle $e^n\p\U$ after time $\tau$.
Conditioned on $\neg\ev U(\hat\alpha)$,
the probability that $\hat \alpha[\sigma,\tau]$ does 
not disconnect $0$ from $\infty$
is bounded by $c_0'e^{-v_0'm/2}$
by (\ref {discon}).
Hence, we get 
$$
1_{\neg\ev U (\hat \alpha)} Z_n
\le c_0' e^{-v_0'm /2 }  
{\ZZ}_n (\eta^n, \hat \beta^n)
.$$
The same applies to $1_{\neg\ev U (\hat \beta)} Z_n$.
Therefore, by  (\ref {up2}),
$$
\EB { 
\bigl(Z_n^\lambda +  Z_n'{}^\lambda \bigr)\,1_{\neg\ev U} }
\le c e^{-v_0'm\lambda/2 } e^{-n \xi}.
$$
It remains to study 
$$
\EB{ \bigl| Z_n^\lambda - Z_n'{}^\lambda \bigr|\,1_{\ev U }  }.
$$

\medskip

\noindent
{\sc Step 2.}
We now show that 
\begin {equation}
\label {ok2}
1_{\ev U}\, \bigl|Z_n - Z_n'\bigr|
\le
c e^{-vm} {\ZZ}_n (\hat \gamma^n)\,,
\end {equation}
assuming $v\le\min\{v_0,v_0'\}$.
When $\ev U$ is satisfied, then the contribution to 
$|Z_n - Z_n'|$ comes from two possible events:
the coupling between $B$ and $B'$ does
not succeed  (this 
occurs with probability at most $ce^{-m v_0}$, independently
from $\hat \gamma$, $\hat B$ and $\hat B'$),
or the coupling succeeds, but
 $\hat B^n$ 
visits
$e^{-m}\U$ and feels the difference between 
$\gamma$ and $\gamma'$.   
In the latter case, the Brownian motion $\hat B$ has 
a probability at most $c_0'e^{-v_0'm}$ not to disconnect the 
origin from infinity after the first visit
to $e^{-m}\p\U$ and before the next visit to
$\p\U$, and after that, it is again an ordinary Brownian
motion up to its hitting time of the circle $e^n\p\U$.
{}From this, (\ref {ok2}) follows.

\medskip

\noindent
{\sc Step 3.}
Suppose first that $\lambda \ge 1$. Then,
since $\max\{ Z_n, Z_n'\} \le {\ZZ}_n (\hat \gamma^n)$,
\begin {align*}
\EB{ 1_{\ev U}\, \bigl|Z_n^\lambda - Z_n'{}^\lambda\bigr|}
&\le 
\lambda\,  \EB{ 1_{\ev U}\, |Z_n - Z_n'|\, {\ZZ}_n^{\lambda-1}}
\\
&\le 
c\,\lambda\, e^{-vm} \E[ {\ZZ}_n^\lambda ]
\qquad\hbox{(by \eref{ok2})}
\\
& \le 
 c'\, \lambda\,  e^{-vm} e^{-n \xi}
\qquad\hbox{(by \eref{up2})}.
\end {align*}
When $\lambda \le 1$, 
\begin {align*}
\EB{ 1_{\ev U} \,\bigl| Z_n^\lambda - Z_n'{}^\lambda\bigr|}
& \le 
\EB{1_{\ev U} \,\bigl|Z_n- Z_n'\bigr|^\lambda }
\\
& \le 
c^\lambda \, \Eb{e^{- \lambda vm}\, {\ZZ}_n^\lambda }
\\
&\le 
c' e^{- \lambda vm} e^{-n \xi }.
\end {align*}
This concludes the proof of Proposition \ref {ok3}
with $v_1= \min \{ v_0'\lambda/2, v_0', v_0, v_0 \lambda\}$. \qed

\subsection {Proof of Proposition \ref {OK}(i)}

Fix $\lambda >0$. 
Define for all integers $k \ge 0$,  
and all $f \in  \A$,  $\gamma \in \Gamma$,
$$
U_k f (\gamma) = \E \left[ f(\bar \gamma )
\frac {\psi^k }{k!} e^{-\lambda \psi} \right].
$$
Taking $y = \lambda \psi$ in the 
inequality $|y|^k / k! \le e^{|y|}$ gives
 for all 
$f \in \A$ and $\gamma \in \Gamma$,
\begin {equation}
\label {U<}
| U_k f( \gamma) |
\le \E [f (\bar \gamma) \lambda^{-k} ]
\le \| f \| \lambda^{-k}.
\end {equation} 
Hence, by  dominated convergence,
for all $z \in \C$ such that $|z| < \lambda$,
for all $f \in \A$ and $\gamma \in \Gamma$,
$$
T_{\lambda - z} f (\gamma) = 
\sum_{k=0}^\infty z^k U_k f(\gamma)
.$$
We will show that there are constants $a,c,v>0$ such
that $N_u(U_k) \leq c a^k$
for all $u\in(0,v)$ and all $k\in\N$.
{}From this, Proposition \ref {OK}(i) immediately follows. 

To find an upper bound for $N_u (U_k)$, first note
that $\| U_k \| \le \lambda^{-k}$ 
by (\ref {U<}) and 
consequently, 
 for all $m < 24$, for all $(\gamma, \gamma') \in \ev X_m$,
$$
\bigl| U_k f(\gamma) - U_k f(\gamma')\bigr|
\le 2 \| f \| \lambda^{-k}.
$$

Suppose that $m \ge 24$ and
$(\gamma , \gamma') \in \ev X_m$. 
As in Proposition \ref{ok3}, since $\gamma$ and $\gamma'$ have the same 
endpoints, we can choose $\hat \gamma = \hat \gamma'$.
Define 
  the event $\ev U'$ that neither $\hat \alpha$
nor $\hat \beta$ hit the circle $e^{- m/24} \p \U$
before $e \p \U$.
Note that when $\ev U'$ is satisfied and $\bar \gamma \in 
\Gamma$, then $(\bar \gamma , \bar \gamma') \in \ev X_m$
(this is where the assumption $m \ge 24$  is needed), so that 
$$
1_{\ev U'}
| f( \bar \gamma) - f (\bar \gamma') | \le \| f \|_u e^{-um}.
$$
An elementary computation 
shows  that if $k \in \N$ and $x \in (0,1)$
\begin{equation}\label{bs}
\frac{1}{k!}\,\Bigl|  \frac d{dx}\bigl({ x^2 (\log x)^k}\bigr)  \Bigr| 
    \le 3,
\end{equation}
so that 
$$
\frac {1}{k!} \,
\EB{
\Bigl | 
\bigl(\frac{\lambda \psi}{2}\bigr)^k e^{-\lambda \psi }
-
\bigl(\frac{\lambda \psi'}{2}\bigr)^k e^{- \lambda \psi'} \Bigr| }
\le 
3 
\EB{ \bigl| e^{-\lambda \psi/2 } - e^{- \lambda \psi' /2} \bigr| }.
$$
 Hence, for all $f \in \A_u$ and for
 all $(\gamma,
\gamma') \in \ev X_m$ (with $m \ge 24$),
\begin {eqnarray*}
\lefteqn {\bigl| U_k f(\gamma) - U_k f(\gamma')\bigr|} \\
& \le 
& \EB{ (1_{\ev U'
} + 1_{\neg \ev U'})
\bigl| f(\bar \gamma) - f(\bar \gamma')\bigr| 
 \frac {\psi'{}^k e^{- \lambda \psi'}}{k!}  }
+
\| f \|\, \frac 1{k!} \,\EB{\bigl|
\psi^k e^{- \lambda \psi}-
\psi'{}^k e^{- \lambda \psi'}
\bigr| } \\
&\le
&\| f \|_u\, e^{-um} \lambda^{-k} +
2 \| f \| \lambda^{-k} \EB { 1_{\neg \ev U'} 
1_{\bar \gamma' \in \Gamma}}
\\
&& \null \qquad +  
\| f \|\:  3 
({2} / {\lambda} )^k 
 \: 
\EB {\bigl| e^{- \lambda \psi/2 } - e^{- \lambda\psi'/2 } \bigr| }
\\
&\le 
& c \| f \|_u ( 2 / \lambda)^k e^{-um}
\end {eqnarray*}
for all  $u \le \min \{ v_1 (
 \lambda/2), v_0'/24\}$,
 by Proposition~\ref{ok3} and (\ref {discon}).

Finally, combining the above estimates shows that 
 for all $u$ smaller than both $v_1 (\lambda/2)$ and  $v_0'/24$,
there exists $c >0$ such that for all $k \ge 0$,
$N_u(U_k) \leq c (2/\lambda)^k$.
This completes the proof
of Proposition \ref {OK}(i).
\qed
 
\section {Spectral gap}  \label{specsec}

We now study the spectrum of $T_\lambda$
for fixed $\lambda > 0$.
The proof of the existence of a spectral gap will be based
on a coupling argument.

\subsection {Coupling the weighted paths}
Let $\r>0$.  For $\hat\gamma^\r$-measurable 
events $\ev A$ define the weighted probability
measures
\begin{equation}\label{Pwdef}
\Pw\gamma/\r/[\ev A]:= \frac{\Eb{ 1_{\ev A}\,
Z^\lambda_\r}}{\Es{Z^\lambda_\r}}\,.
\end{equation}
For $0<k\le n$, 
let $\delta(n,k,\gamma)$ denote a random variable with the
same law that $\bg^k$  has under  $\Pw\gamma/n/$.
It is easy to verify that
\begin{equation}\label{convol}
\delta\bigl(n-j,k,\delta(n,j,\gamma)\bigr)
\hbox{ has the same law as }
\delta(n,k+j,\gamma)\,,
\end{equation}
when $j+k\le n$.

The following coupling result will be crucial in our proof:

\begin {proposition}
\label {coupling}
There exist constants $v_2,c>0$ such that for all $n \ge 1$,
for all $\gamma,\gamma'\in\Gamma$, 
one can define $\delta=\delta(n,n,\gamma)$ and $\delta'=\delta(n,n,\gp)$ 
on the same probability space $(\Omega,
{\cal F},\coupl)$ such that 
$$
\Bb\coupl{ (\delta, \delta') \notin \ev X_{n/3} }
\le  c\, e^{-v_2n}.
$$
\end {proposition}

This subsection will be devoted to the proof of 
this proposition.
The rough strategy of the proof is as follows.  The first
step is to get both paths to be in $\Gamma^+$ (recall
the definition of $\Gamma^+$ from Section~\ref{Ssep}).
  The second step
is to make the two paths match up and walk together
a little while.  The third step is to show that if
the two paths have walked together for some time,
then it is unlikely they will decouple.  If any of
these steps fails, the coupling process returns to
the beginning.  (For technical reasons, the order
in which we address these steps is different from the logical
order indicated above.)

Suppose that $u>0$ is fixed (and small enough), and choose
$C$ such that for all $\gamma \in \Gamma^+$,
$ R_1 (\gamma) \ge C $.
Define for all $k \ge 1$,  
$$
\K_k 
= \bigl\{ ( \gamma, \gamma') \in \Gamma^2 
\st R_1 (\gamma) \ge Ce^{- ku}, \ R_1 (\gamma') \ge Ce^{-ku}
\hbox { and } (\gamma , \gamma') \in \ev X_k \bigr\}.
$$

\begin {lemma}
\label {coupling2}
There exist $c,w>0$ such that for all $n \ge 1$,
for all $k \ge 24$, for all $(\gamma, \gamma') \in \K_k$, 
one can define $\delta=\delta(n,1,\gamma)$ and $\delta'=\delta(n,1,\gp)$ 
on the same probability space $(\Omega, {\cal F}, \coupl)$
such that 
$$
\Bs\coupl{  (\delta, \delta') \in \K_{k+1} } \ge 1- ce^{-wk}.
$$
\end {lemma}

Let us recall the following 
straightforward fact (see, e.g., \cite {Li}). 
Suppose that $\nu$ is a probability measure on a 
measurable space $(S, {\cal S})$, and that $\nu_1$ and $\nu_2$ are
two other probability measures on this space such that the 
Radon-Nikodym derivative $q(\cdot)$ (resp.\ $q' (\cdot)$) 
of $\nu_1$ (resp.\ $\nu_2$) 
with respect to $\nu$ exists.
Then, one can find a probability space
$(\Omega , {\cal F}, \coupl)$ and measurable functions
$\delta: \Omega \to S$ and $\delta' : \Omega \to S$, such that the
law of $\delta$ is $\nu_1$, the law of $\delta'$ is $\nu_2$, and
$$ \Bs\coupl{  \delta = \delta' } \ge \int_S \min (q,q') \: d\nu 
= 1 - \frac {1}{2} \int | q - q'|  \,d \nu = 1-\frac12 \|\nu_1-\nu_2\| .$$
A similar argument was implicitly used in the proof of 
 \ref{lemma.apr26}.(ii).

\proofof{Lemma \ref {coupling2}}
Suppose that $(\gamma, \gamma') \in \K_k$ with $k \ge 24$. 
The law of $\delta$ has Radon-Nikodym derivative
$$
q:=
 \frac { e^{\xi-\lambda \psi_1} R_{n-1} (\bg)}{R_n (\gamma)}
$$
with respect to the law of $\bg$
and similarly for $\delta'$:
$$
q':=
 \frac { e^{\xi-\lambda \psi_1'} R_{n-1} (\bgp)}{R_n (\gp)}\,.
$$
Since $\hat\gamma^1$ has the same law as
${\hat\gamma}'{}^1$, we may take them
to be the same.
Our first  goal is to estimate $\Es {\left|q-q'\right|}$.

Recall that $R_{n-1}( \bar \gamma') \le c_1$ by (\ref {up2}), and that 
$R_n (\gamma') \ge (C/c_0) e^{-uk}$
by (\ref{use})  and the fact that
 $(\gamma,\gamma') \in \K_k$.
Hence, by Proposition~\ref{ok3},
$$
\EB{ \frac { R_{n-1} (\bar \gamma')}
{R_n (\gamma')}
\bigl| e^{- \lambda \psi } - e^{- \lambda \psi'}\bigr|
}
\le c e^{uk}  e^{-v_1 k } \le c e^{-uk}
$$
when $u \le v_1/2$.
On the other hand,
using Proposition \ref {ok3} again, it follows that 
\begin {align*}
\EB{ 1_{(\bar \gamma, \bar \gamma') \in \ev X_k}
  \frac {\bigl| R_{n-1} (\bar \gamma') - 
R_{n-1} (\bar \gamma)\bigr| }{ R_n (\gamma')}
\,e^{- \lambda \psi}}
& \le 
c e^{uk}\, \Eb{1_{(\bar \gamma, \bar \gamma') 
\in \ev X_k}
  | R_{n-1} (\bar \gamma) - R_{n-1} (\bar \gamma')|}
\\&
\le c' e^{uk} e^{- v_1 k}
\le c' e^{-uk}
\end {align*}
when $u \le v_1/2$.
Since $k \ge 24$ and $(\gamma, \gamma') \in 
\ev X_k$, if $(\bar \gamma , \bar  \gamma') \in \Gamma\times\Gamma\setminus
\ev X_{k+1}$, then $\hat \gamma^1$ has a 
downcrossing from $1$ to $e^{-k/24}$.
This shows readily that
$$
\EB{
1_{(\bar \gamma, \bar \gamma') \notin \ev X_k}
\frac {R_{n-1} (\bar \gamma') + R_{n-1} (\bar \gamma)}
{R_n (\gamma')} }
\le  c e^{u k} e^{-k v_0' / 24}
\le c e^{-k u}$$
for all $u < v_0' / 48$.
Finally, 
\begin {align*}
\EB{ 
R_{n-1} (\bar \gamma) 
 \bigl| {R_n (\gamma)^{-1}} - {R_n (\gamma')^{-1}} \bigr|
\, e^{- \lambda \psi}
}
& \le 
c e^{2uk} \bigl|R_n (\gamma) - R_n (\gamma')\bigr|
\\
& \le 
c' e^{2uk}  e^{-kv_1}
\le  c' e^{-uk}
\end {align*}
for all $u \le v_1/3$.

Putting the pieces together, we see that if we take $u<
\min ( v_0' / 48,  v_1/3)$
then $\Eb{\left| q-q'\right|}\le c e^{-uk}$,
and hence there is a coupling of $\delta$ and $\delta'$ with
$\Bs\coupl{\delta_{k+1}=\delta'_{k+1}}\ge 1-c e^{-uk}$.
 
We now  check that 
 $
\mu [ R_1 (\delta) \le  C e^{-u (k+1)}]
=\Bb{\Pw\gamma/n/}{R_1 (\bg) \le Ce^{-u (k+1)}}$ is
also  exponentially small in $k$.
Recall that $R_{n-1} (\bar \gamma) \le c_0 R_1 (\bar \gamma)$,
by \eref{use}, so
that 
\begin {align*}
\Bb{\Pw\gamma/n/}{R_1 (\bg) < Ce^{-u (k+1)}}
&
= \EB{
1_{R_1 ( \bar \gamma ) \le C e^{-u (k+1)}  }
\frac { e^{\xi-\lambda \psi} R_{n-1} (\bar \gamma)}{R_n (\gamma)} }
\\ &
\le \EB{ \frac { e^{\xi-\lambda \psi}}{R_n (\gamma)} c_0 C e^{- u (k+1)}}
\\
& = 
c_0 C e^{-u (k+1)} \frac { R_1 (\gamma)}{R_n(\gamma)}
\\ &
\le 
c_0^2 C e^{-u k },
\end {align*}
and a similar inequality holds for $\gamma'$.

Finally, it remains to bound the probability $\mu [
\delta \notin \ev Y_{k+1} ]$.
Since $\gamma \in \ev Y_k$ and $k \ge 24$, if $\delta \notin
\ev Y_{k+1}$ then 
$\hat \gamma^1$ has a downcrossing from 
$1$ to $e^{-k/24}$.
Hence,
\begin {eqnarray*}
\mu [ \delta \notin \ev Y_{k+1} ]
& = &
\EB {
1_{ \bar \gamma \notin \ev Y_{k+1}}
e^{\xi - \lambda \psi} \frac {R_{n-1} (\bar \gamma)}{R_n (\gamma)}}
\\
& \le &
2c_0' e^{- v_0'k/24}  (c_0 c_1 /C) e^{uk} e^{\xi} \\
& \le &
ce^{-uk}
\end {eqnarray*}
when $u < v_0' / 48$.
This completes the proof of the lemma. \qed

\medbreak
We now choose $p \ge 25$ such that 
(for the constants defined in Lemma \ref {coupling2}), 
$ ce^{-w(p-1)}<1/2$.
This is to make sure that the coupling in Lemma \ref {coupling2}
occurs with positive probability for all $k \ge p -1$.

\begin {lemma}
\label {coupling1}
There exists a constant $c=c(p)>0$
such that for all $\gamma, \gamma' \in \Gamma^+$,
and for all $n \ge p$,
one can define $\delta=\delta(n,p,\gamma)$ and $\delta'=\delta(n,p,\gamma')$
on the same probability space $(\Omega, {\cal F}, \coupl)$, such that 
$$
\Bb\coupl{  (\delta, \delta') \in \K_{p-1} } \ge c .
$$
\end {lemma}

\proof
Take $\gamma,\gp\in\Gamma^+$.
Define the Brownian motions $\hat \alpha$ and $\hat \alpha'$
on the same probability space by mirror coupling.
That is, we take $|\hat \alpha(t)|=|\hat \alpha'(t)|$ 
and keep $\arg\hat\alpha(t)+\arg\hat\alpha'(t)$ constant
up to the first time $t$ at which $\hat\alpha(t)=\hat\alpha'(t)$.
After they have met, they stay together.
Couple  $\hat \beta$ and $\hat \beta'$ in the same way.

It is easy to see that there is a $c=c(p)>0$,
which does not depend on $\gamma$ or $\gp$, such that 
with probability at least $c$,
\begin {itemize}
\item
$\hat \gamma$ and $\hat \gamma'$ coalesce before they reach $e\p\U$,
\item
$\bg^p,\bgp^p\in\Gamma^+$,
\item 
$(\bg^p, \bgp^p) \in \ev X_{p-1}$
and
\item 
$\hat\gamma^p\cup\hat \gamma'{}^p 
\subset
\bigl\{\, r e^{i \theta} \st   r > e^{-1/8}, \, \theta  \in ( 
\pi/4, 7\pi/4 ) \bigr\}.
$
\end {itemize}  
Let $\ev G$ denote this event.
The law of $\delta$ has Radon-Nikodym derivative
$$
q:=
 \frac { e^{p\xi-\lambda \psi_p} R_{n-p} (\bg^p)}{R_n (\gamma)}
$$
with respect to the law of $\bg$
and similarly for $\delta'$:
$$
q':=
 \frac { e^{p\xi-\lambda \psi_p'} R_{n-p} (\bgp^p)}{R_n (\gp)}\,.
$$
Since $\P[\ev G]$ is bounded from below,
to prove that there is a constant $c>0$ such that
for all $n\ge p$ and all $\gamma,\bg\in\Gamma^+$ there
is a coupling $\coupl$ with 
$$
\Bb\coupl{\{\delta_{p-1}=\delta_{p-1}'\}
\cap\{\delta,\delta'\in\Gamma^+\}}\ge c
,
$$
it suffices to show that $q$ and $q'$ are bounded away
from $0$ on $\ev G$.  This does hold, since $R_n$ is bounded
and bounded from zero on $\Gamma^+$, and
one can verify directly that $\psi_p$ and $\psi'_p$ are
bounded on $\ev G$.
Because $R_1\ge C$ on $\Gamma^+$, such
$(\delta,\delta')$
are in $\K_{p-1}$, and hence the proof is now complete.
\qed

\proofof{Proposition \ref {coupling}}
Set $p^*=p+1$.
Suppose first that $n=mp^*$, where $m\in\N$.
Let $\gamma, \gamma'\in\Gamma$.
The coupling $\coupl$ is defined as follows.
Inductively, we construct a sequence $\bigl(\delta({j}),
\delta'({j})\bigr)$,
$j=0,1,\dots,m$, such that 
$\delta(0):=\gamma$, $\delta'(0):=\gamma'$,
$\delta({j+1})$ has the law of $\delta\bigl(n-jp^*,p^*,\delta(j)\bigr)$ and
$\delta'(j+1)$ has the law of $\delta\bigl(n-jp^*,p^*,\delta'(j)\bigr)$.
Then, we set $(\delta,\delta'):=\bigl(\delta(m),\delta'(m)\bigr)$. 
Repeated use of \eref{convol} shows that $\delta$ and
$\delta'$ have the desired laws.

Let 
$$K_{j} 
:= \max \Bigl\{ k \ge 1 \st \bigl(\delta({j}), \delta'({j})\bigr) \in \K_k
\Bigr\},$$
and let $K_j:=0$ if the set on the right hand side is empty.

It follows easily from Lemma \ref {coupling2} iterated $p^*$
times that  it is possible to 
construct $\bigl(\delta({j+1}), \delta'({j+1})\bigr)$ 
in such a way that 
$$
\PB{  K_{j+1} \ge  K_j + p^* \md  \delta(j),\delta'(j)}
\ge \bigl( 1- c e^{-wK_j}\bigr)^{p*}\,1_{K_j\ge p-1}\,.
$$
For the case where $K_j< p-1$, the construction of
$\bigl(\delta(j+1),\delta'(j+1)\bigr)$ proceeds as follows. 
Note that the strong Markov property, the Separation Lemma~\eref{sep},
and \eref {tildec} imply that
$$
\inf_{\r\ge 1}\inf_{\gamma\in\Gamma}\Pw\gamma/\r/[\bg\in\Gamma^+]>0\,.
$$
Therefore, Lemma \ref {coupling1}
show that it is 
possible to construct $\bigl(\delta({j+1}), \delta'({j+1})\bigr)$ in
such a way that  
$$
\inf_{n,j,\gamma,\gamma'}\
\inf \PB { K_{j+1} \ge p-1 \md \delta(j),\delta'(j)}>0\,.
$$
 
By comparison with a Markov chain on the integers (see, e.g.,
Proposition 2.(iii) in \cite {BFG}), this implies
readily that 
$$
\Pb{  K_{m} \ge p^*m/2 } \ge 1 - c' e^{-mw'}
$$
for some constants $c',w' >0$.
This proves the Proposition when $n/p^*\in\N$.
The general case follows easily.
For instance, if
 $n = m p^* + m'$ with 
$m' \in [0, \cdots, p^*)$, we can apply the result 
for $n' = m p^*$ to $\delta= \delta (n, m', \gamma)$ 
and to $\delta' = \delta ( n ,m', \gamma')$, and note that 
$ \delta( n', n', \delta ) 
= \delta (n, n, \gamma)$ and $\delta (n',n', \delta')
= \delta  (n,n,\gamma')$.
The small values of $n$ 
can be handled by modifying the constant $c$ in the statement
of the lemma.
\qed

\subsection {Proof of Proposition \ref {OK}.(ii)}

We now conclude the proof of Proposition \ref {OK}(ii)
and thereby the proof of Theorem~\ref {main}.

\noindent
{\sc Step 1.} 
Let $k>0$.
Suppose that $\gamma \in \Gamma$ is fixed,
and let $\gamma':=\delta(n+k,k,\gamma)$.
Set $\delta:=\delta(n,n,\gamma')$
and $\delta':=\delta(n,n,\gamma)$.
Then $\delta$ has the law of $\delta(n+k,n+k,\gamma)$.
By Proposition~\ref{coupling}, $\delta$ and
$\delta'$ may be defined on the same
probability space $(\Omega, {\cal F}, \coupl)$, so that
$$
\Bb\coupl { (\delta,  \delta') \notin \ev X_{n/3} } 
\le c e^{-nv_2}.
$$
Hence, for all $f \in \A_u$,
\begin{equation}
\label{au}
\begin{aligned}
\left|
\frac { T^n f(\gamma) }{T^n 1 (\gamma)} - 
\frac {T^{n+k} f(\gamma)} {T^{n+k} 1 (\gamma)}
\right|
& \le  
\int| f (\delta) - f (\delta') | \,d\coupl
\le 
 \| f \|_u ( 2ce^{-nv_2} + e^{-nu/3} ) 
. \end {aligned}
\end{equation}
It follows that for all $f \in \A_u$, $T^n f(\gamma) / T^n1(\gamma)$
converges when $n \to \infty$ to some limit $h(f, \gamma)$.

The same kind of argument gives for all $\gamma,\gamma'\in\Gamma$,
$$
\left|
\frac { T^n f(\gamma) }{T^n 1 (\gamma)}
-
\frac {T^{n} f(\gamma')} {T^{n} 1 (\gamma')}
\right|
 \le 
\| f \|_u ( 2ce^{-nv_2} + e^{-nu/3} ) ,
$$
and therefore
the limit $h(f, \gamma) = h(f)$ is in fact independent of $\gamma$.
Clearly, $h : \A \to \C$ is linear and 
$|h(f)| \le \| f \|$ for all $f 
\in \A$ so that $h$ is a  bounded linear functional on 
$\A_u$.
 
\medbreak
\noindent
{\sc Step 2.}
We are going to find an upper bound for  the operator
norm $N_u $
of the operator $f\mapsto T^nf-h(f)T^n1$.
Inequality~\eref{au} shows that for all $f \in \A_u$ and $\gamma \in
\Gamma$,
$$
\bigl| T^nf(\gamma) - h(f) T^n1 (\gamma)
\bigr| \le  \| f \|_u (2c e^{-nv_2}+ e^{-nu/3}) T^n1 (\gamma)
\le c' \| f\|_u e^{-n ( \xi + u/3)}
$$
for all sufficiently small $u \le 3 v_2$.
Suppose now that $(\gamma , \gamma') \in \ev X_m$ 
and that $m \le n/6$.
Then, the previous estimate gives
$$
\bigl| T^n f(\gamma) - h(f) T^n (\gamma) -
T^n f (\gamma') + h(f) T^n (\gamma') \bigr|
\le 2 c' \| f \|_u
e^{-n \xi} e^{- nu /6} e^{-mu}
.$$
Assume now that $(\gamma, \gamma') \in \ev X_m$ and that 
$m \ge n/6$.
Defining $\hat \gamma = \hat \gamma'$, and using Proposition~\ref{ok3}
and the fact that $|h(f)| \le \| f \|$,
we get 
\begin {eqnarray*}
\lefteqn
{ \bigl|T^n f (\gamma) - T^n f(\gamma')\bigr| +
\bigl| h(f)\bigr| \bigl|T^n1 (\gamma) - T^n 1 (\gamma')\bigr| }
\\ &
\le & 
\Eb{ \bigl| f(\bar \gamma^n)e^{- \lambda \psi_n} -
  f(\bar \gamma'{}^n)e^{- \lambda \psi'_n} \bigr| }
+   \| f \| \,\Eb {\bigl| e^{- \lambda \psi_n'} - e^{-\lambda \psi_n}
\bigr|}
\\ &
\le & 
\Eb{ \bigl| f(\bar \gamma^n) - f(\bar \gamma'{}^n) \bigr| 
e^{- \lambda \psi_n} }
+  2 \| f \| \Eb {\bigl| e^{- \lambda \psi_n'} - e^{-\lambda \psi_n}
\bigr|}
\\ &
\le &  
\EB {
1_{\ev V} \bigl| f(\bar \gamma^n) - f(\bar \gamma'{}^n) \bigr| 
e^{- \lambda \psi_n} }
+
2 \| f \| \EB {
1_{\neg \ev V} 
e^{- \lambda \psi_n} }
+
  2 \| f \| \Eb {\bigl| e^{- \lambda \psi_n'} - e^{-\lambda \psi_n}
\bigr|} 
,
\end {eqnarray*}
where $\ev V$ is the event that $(\bar \gamma^n, \bar \gamma'{}^n)
\in \ev X_{n+m}$.
We are going to bound the three terms separately.
The last one is bounded by 
$ 2c \| f \| e^{-n \xi} e^{-m v_1}$ 
(and therefore by $ 2 c \| f \| e^{-n \xi}
e^{-mu} e^{-nu/6}$ for all $u < v_1/2$)
by 
Proposition~\ref {ok3}. 
For the first one:
\begin {eqnarray*}
\EB {
1_{\ev V} \bigl| f(\bar \gamma^n) - f(\bar \gamma'{}^n) \bigr| 
e^{- \lambda \psi_n} }
&\le & \| f\|_u e^{-(n+m)u} \EB { e^{-\lambda \psi_n}} 
\\
& \le&  c' \| f\|_u e^{-n \xi} e^{-mu} e^{-nu}.
\end {eqnarray*}
For the second term, note that on
$\bigl\{\bg^n\in\Gamma\bigr\}\setminus \ev V$
there is a
$j\in\{0,1,2,\dots,23\}$ such that
$\bg^n$ has a downcrossing from  
$e^{-jn/24}$ to $e^{-(jn+m)/24}$.  
By estimating the contribution of each of these 24 possible
values of $j$
separately, one easily  gets, using the strong Markov property,
(\ref {discon}) and (\ref {up2}) 
  that 
\begin {eqnarray*}
 \EB {
1_{\neg \ev V} \,
e^{- \lambda \psi_n} }
&\le&
\sum_{j=0}^{23}
\EB { e^{- \lambda \psi_{jn/24}}\, 2c_0'\, e^{-v_0'm / 24}\,
c_1 \,e^{-(24-j)n \xi / 24} }\\
& \le &
  c\, e^{-n\xi} \,e^{-v_0'm  /24}
\\
& \le &
c \,e^{-n \xi} \,e^{-u m} \,e^{-un/6}
\end {eqnarray*}
for all $u < v_0' / 48 $.
Combining these three estimates shows that
for all sufficiently small $u$,
there exists $c =c (u,\lambda)$ such that  
for all $n \ge 1$,
\begin {equation}
\label {graal}
N_u \bigl( T^n (\cdot) - h(\cdot) T^n1 \bigr) 
\le c\, e^{-n \xi } \,e^{-n u /6}
.\end {equation}

\medbreak
\noindent
{\sc Step 3.}
Note that 
$T^{n+1}1 (\gamma) / T^n 1 (\gamma)\to h(T1)$ 
and $T^{n+j} 1(\gamma) / T^n 1 ( \gamma) \to 
h(T^j 1)$ as $n\to\infty$.
Since $T^n 1 (\gamma) \approx e^{-n \xi}$
for $\gamma \in \Gamma^+$,
we get $h( T1) = e^{-\xi}$ and $h(T^j1) =e^{-j \xi}$.
Recall also that $\| T^j1 \|_u \le ce^{-j \xi}$ 
by Proposition \ref {ok3}.
Hence, (\ref {graal}) for $f = T^j1$ shows that 
$$
\| T^{n+j} 1 - e^{-j \xi} T^n 1 \|_u 
\le c\,e^{-n \xi}\, e^{-j \xi }\, e^{-n u/6}.
$$
Hence, $R_n =  e^{n \xi} T^n 1$ converges in $\A_u$ to some limit $R$,
and 
$$
\|R_n-R\|_u = \| e^{n \xi} T^n1 - R \|_u \le c\,e^{-n u/6}.
$$

Since $T$ is continuous, $T\,R=e^{-\xi}R$.
If $h(f)=0$, then $T^nf/T^n1\to 0$ on $\Gamma$
and therefore $T^{n+1}f/T^{n}1\to0$,
which implies $h(Tf)=0$. That is,
$T\,\mathrm{ker}(h)\subset\mathrm{ker}(h)$.
Moreover,  \eref {graal} shows that 
the operator norm of the restriction of $T^n$ to $\mathrm{ker}(h)$ 
is bounded by $c\,e^{-n\xi}\,e^{-nu/6}$;
this implies Proposition~\ref {OK}.(ii). \qed

\ifhyper\def\email#1{\href{mailto:#1}{\texttt{#1}}}\else
\def\email#1{\texttt{#1}}\fi

\bigskip
\filbreak 

\vtop{
\hsize=2.3in
\noindent Greg Lawler\\
Department of Mathematics\\
Box 90320\\
Duke University\\
Durham NC 27708-0320, USA\\
\email{jose@math.duke.edu}
}
\bigskip

\vtop{
\hsize=2.3in
\noindent Oded Schramm\\
Microsoft Corporation,\\
One Microsoft Way,\\
Redmond, WA 98052; USA\\
\email{schramm@microsoft.com}
}
\bigskip

\vtop{
\hsize=2.6in
\noindent Wendelin Werner\\
D\'epartement de Math\'ematiques\\
B\^at. 425\\
Universit\'e Paris-Sud\\
91405 ORSAY cedex, France\\
\email{wendelin.werner@math.u-psud.fr}
}

\filbreak

\end {document}